\documentclass[12pt]{article}

\usepackage{exptex}
\usepackage[top=1.5cm,left=2.5cm,right=2.5cm,bottom=2.0cm,nohead]{geometry}
\DeclareMathAlphabet{\mathitbf}{OML}{cmm}{b}{it}

\begin{document}
\title{Congruence and Metrical Invariants of Zonotopes}
\author{Eugene Gover}
\date{}
\maketitle

The defining matrix $A$ of a zonotope $\mathcal{Z}(A)\subset\mathbb{R}^n$
determines the zonotope as both the linear image of a cube
and the Minkowski sum of line segments specified by the columns of the matrix. A zonotope is also a
convex polytope with centrally symmetric faces in all dimensions.
When a zonotope is represented by a matrix, its volume is the sum of the absolute values
of the maximal-rank minors. Sub-maximal rank minors compute the lower-dimensional volumes of facets.
Maximal-rank submatrices determine various tilings of a zonotope, while those of submaximal rank
define the angles between facets, normal vectors to facets, and can be used to demonstrate
rigidity and uniqueness of a zonotope given various facet-volume and normal-vector data.
Some of these properties are known. Others, are new. They will all be presented
using defining matrices.

The first section focuses on the central symmetry of faces and facets of convex polytopes,
and gives new proofs of theorems of Shephard and McMullen.
The second section introduces the Gram matrix $A^{T}A$, called the shape matrix of the zonotope,
and gives it the central role in a discussion of congruences between zonotopes.
The same matrix also plays an important part in the third section
where new proofs of theorems of Minkowski
and Cauchy-Alexandrov are given in the case of zonotopes.

\smallskip

\section*{1. \,\,\,\,Central Symmetry and Zonotopes}

\bigskip

Given $\boldsymbol{x}$ and $\boldsymbol{c}\in \mathbb{R}^n$, the points $\boldsymbol{x}$ and
$2\boldsymbol{c}-\boldsymbol{x}$ will be said to be {\bf symmetric images} of each other
with respect to $\boldsymbol{c}$.
For a nonempty subset $X\subset\mathbb{R}^n$, the set
$$X_{\boldsymbol{c}}\,\,{\mathrel{\mathop :}=}\,\,2\boldsymbol{c}-X=\left\{2\boldsymbol{c}-\boldsymbol{x}\,|\,\boldsymbol{x}\in X\right\}$$
will be called the {\bf symmetric image} (or {\bf point reflection}) of $X$ with respect to $\boldsymbol{c}$.

$X$ will be called {\bf centrally symmetric} with {\bf center of symmetry} $\boldsymbol{c}$ if and only if there exists
$\boldsymbol{c}\in \mathbb{R}^n$ such that $X=X_{\boldsymbol{c}}$\,, that is, iff $X$ contains the symmetric image
of each of its points with respect to a single center $\boldsymbol{c}$. The condition can be restated as saying
there exists $\boldsymbol{c}\in \mathbb{R}^n$ such that $X=2\boldsymbol{c}-X$, or such that $X-\boldsymbol{c}=-X+\boldsymbol{c}$,
or as the assertion that $\boldsymbol{x}\in X$ iff $2\boldsymbol{c}-\boldsymbol{x}\in X$ holds.
The center of a bounded centrally symmetric set is unique
but need not belong to the set. (If the set is convex, the center will belong to the set.)
Another equivalent condition is that $X$ is centrally symmetric
if and only if there exists a translation $\tau$ such that $\tau (X)=-X$.
The center of symmetry will then be $\boldsymbol{c}=\frac{1}{2}\tau^{-1}(\boldsymbol{0})$. In terms of
$\boldsymbol{c}$, $\tau(\boldsymbol{x})=\boldsymbol{x}-2\boldsymbol{c}=-(2\boldsymbol{c}-\boldsymbol{x})$.
Central symmetry can also be described using a {\bf symmetric cone} over $X$ centered at $\boldsymbol{c}$,
which is defined as the set
$$\text{cone}_{\boldsymbol{c}}\,X\,{\mathrel{\mathop :}=}
\left\{t\boldsymbol{c}+(1-t)\boldsymbol{x}\,|\,\boldsymbol{x}\in
X,\,0\leq t\leq 2\right\}.$$
The subset of cone$_{\boldsymbol{c}}\,X$ for a single value $t_0\in [0,2]$ will be denoted
$$\text{cone}_{\boldsymbol{c},{\scriptstyle t_0}}\,X{\mathrel{\mathop :}=}
\left\{t_0\boldsymbol{c}+(1-t_0)\boldsymbol{x}\,|\,\boldsymbol{x}\in
X\right\}.$$
In particular, \,cone$_{\boldsymbol{c},0}\,X\,=\,X$,
\,cone$_{\boldsymbol{c},1}\,X=\{\boldsymbol{c}\}$, \,and
\,cone$_{\boldsymbol{c},2}\,X=\left\{2\boldsymbol{c}-\boldsymbol{x}\,|\,\boldsymbol{x}\in X\right\}
=X_{\boldsymbol{c}}$.
$X$ will then be centrally symmetric with respect to $\boldsymbol{c}$ iff
cone$_{\boldsymbol{c},0}\,X=\text{cone}_{\boldsymbol{c},2}\,X$.

The following properties are easily verified:

\medskip

\noindent {\bf Lemma 1.1.}  {\it
$(a)\hspace*{3.6bp}(X_{\boldsymbol{c}})_{\boldsymbol{c}}=X$;$\,\,\,(X_{%
\boldsymbol{c}_1})_{\boldsymbol{c}_2}=2(\boldsymbol{c}_2-%
\boldsymbol{c}_1)+X$;$\,\,\,((X_{\boldsymbol{c}_1})_{\boldsymbol{c}_2})_{%
\boldsymbol{c}_3}=X_{\boldsymbol{c}_3-\boldsymbol{c}_2+\boldsymbol{c}_1}$;

$(((X_{\boldsymbol{c}_1})_{\boldsymbol{c}_2})_{\boldsymbol{c}_3})_{%
\boldsymbol{c}_4}=2(\boldsymbol{c}_4-\boldsymbol{c}_3+\boldsymbol{c}_2-%
\boldsymbol{c}_1)+X$, etc.

$(b)$\hspace{4.6bp}For any $\boldsymbol{c}\,\in\mathbb{R}^n$,} cone$_{\boldsymbol{c}}\,X$
{\it is centrally symmetric with center of symmetry $\boldsymbol{c}$.

$(c)$\hspace{4.6bp}For any $\boldsymbol{c}\,\in\mathbb{R}^n$, $X\cup X_{\boldsymbol{c}}$ is
centrally symmetric with center of symmetry $\boldsymbol{c}$.

$(d)\hspace{4.6bp}X$ is centrally symmetric iff $X_{\boldsymbol{c}}$ is centrally symmetric
for every $\boldsymbol{c}$.

$(e)\hspace{3.6bp}X$ is centrally symmetric iff for each $\boldsymbol{c}$ there exists
$\boldsymbol{v}\,$ such that $X_{\boldsymbol{c}}=\,\boldsymbol{v}+X$.

$(f)\hspace{2.6bp}$If $X$ and $Y$  are each centrally symmetric with respect to the same center $\boldsymbol{c}$,

then $X\cap Y$ and $X\cup Y$ are also centrally symmetric with respect to $\boldsymbol{c}$.}

\medskip

Part$\,(a)$ says two successive point reflections with the same center leave a set
unchanged while using different centers results in translation by twice the
difference between the centers. More generally, an odd number of point reflections
with centers $\boldsymbol{c}_1,\ldots,\boldsymbol{c}_{2n+1}$ are equivalent to
a single reflection with respect to the alternating sum of the centers; in
particular, the image of $\boldsymbol{x}\in X$ will be
2$\left(\boldsymbol{c}_{2n+1}-\cdots+\boldsymbol{c}_1\right)-\boldsymbol{x}%
\in X_{(\boldsymbol{c}_{2n+1}\,-\,\cdots\,+\,\boldsymbol{c}_1)}$. An even
number of reflections are equivalent to translation by twice the alternating
sum of the centers. Parts $(b)$ and $(c)$ say that the symmetric cone of a
set and the union of the set with any point reflection are centrally symmetric.
Part $(d)$ says that a set is centrally symmetric if and only if every point
reflection is centrally symmetric. Part $(e)$ says $X$ is centrally symmetric
iff it can be translated to any and every point reflection image of itself.

\bigskip

Consider a unit cube positioned along the coordinate axes of ${\mathbb R}^k$. Its image in ${\mathbb R}^n$
under a linear transformation defined with respect to the standard bases by a real $n\times k$ matrix $A$ with columns
$\boldsymbol{a}_1,\ldots,\boldsymbol{a}_k$ is the set
$\mathcal{Z}(A)\,=\,\mathcal {Z}({\mathitbf a}_1,\ldots, {\mathitbf a} _k)\,=\,\{\sum t_i{\mathitbf a}_i\,|\,0\leq t_i\leq 1\}$.
The set will be called the {\bf zonotope} generated by the columns of $A$, which will in turn be called
the {\bf defining matrix} of $\mathcal{Z}(A)$. The {\bf rank} of the zonotope is the rank of its defining matrix.
In the special case when $n\geq k$ and the columns are independent ({\it i.e.}, rank $A=k$), the image is also a {\bf parallelotope} and can be denoted as $\mathcal{P}(A)$
or $\mathcal{P}({\mathitbf a}_1,\ldots, {\mathitbf a} _k)$.
The parallelotopes we will consider are generated by the independent columns of tall and thin, or square matrices of rank $k$ with $n\geq k$.
They are skewed, stretched, or shrunken images of cubes.
The zonotopes that are not parallelotopes will be generated by the dependent columns of matrices
of rank $r$ with $r < k$. They are flattened images of cubes.

The {\bf Minkowski sum} of
sets of $S_1,\ldots,\,S_k\subset\mathbb{R}^n$, denoted with the symbol $\oplus$, is the set
$S_1{\oplus}\cdots {\oplus}S_k\,=\,\{{\mathitbf s}_1+\cdots+{\mathitbf s}_k\,|\,{\mathitbf s}_i\in S_i\}$.
A zonotope is the Minkowski sum of line segments:
$\mathcal{Z}(A)\,=\,l\boldsymbol{a}_1{\oplus}\cdots{\oplus}l\boldsymbol{a}_k$
where the line segment $l\boldsymbol{a}_i\,=\,\{t\boldsymbol{a}_i\,|\,0\leq t\leq1\}$.
For a parallelotope, the generators are linearly independent and the Minkowski sum of the corresponding line segments yields a prism whose
base is any Minkowski sum leaving out one of the segments. (Note that the parallelotopes we will consider form a proper subset
of the polytopes that fill space by translation, which are also called parallelotopes.)
Cubes are convex, centrally symmetric, and the convex hulls of finite point sets. It follows that zonotopes, which are their linear images,
are also convex, centrally symmetric polytopes. (As polytopes, zonotopes are also finite intersections of half-spaces.)
Regarded as Minkowski sums of line segments, zonotopes are centrally symmetric for another reason: for each
$\sum_i t_i\boldsymbol{a}_i\in\mathcal{Z}(A)$, there corresponds
$\sum_i (1-t_i)\boldsymbol{a}_i\in\mathcal{Z}(A)$; these two points
have center of symmetry $\sum_i\boldsymbol{a}_i /2$, which becomes the center of the entire zonotope.

Suppose zonotope $\mathcal{Z}(A)$  is defined by the matrix
$[\boldsymbol{a}_1,\ldots ,\boldsymbol{a}_k]\in\mathbb{R}^{n\times k}$ of rank $r\leq k$.
A  subzonotope of rank $s\leq r$ of the form
$\mathcal{Z}(\boldsymbol{a}_{j_1},\ldots,\boldsymbol{a}_{j_t})$ to which no further column vectors
can be added as generators without increasing the rank
will be called a {\bf generating face of dimension $s$} or an {\bf$s$-face} of $\mathcal{Z}(A)$.
The generators themselves are considered the {\bf generating 0-faces}.
A line segment $\mathcal{Z}(\boldsymbol{a}_i)\,=\,l\boldsymbol{a}_i$
will be a generating 1-face or {\bf edge} unless there is a larger, maximal collection
${\boldsymbol{a}_{i_1},\ldots ,\boldsymbol{a}_{i_t}}$ of generators containing $\boldsymbol{a}_i$
with each generator a scalar multiple of the others. In that case,
$\mathcal{Z}(\boldsymbol{a}_{i_1},\ldots ,\boldsymbol{a}_{i_t})$ becomes a generating edge of $\mathcal{Z}(A)$
containing each of the $\mathcal{Z}(\boldsymbol{a}_i)$'s.
A generating $(r-1)$-face will be called a {\bf generating facet}.
A rank $r$ subzonotope with exactly $r$ generators will be called a
{\bf generating parallelotope} of $\mathcal{Z}(A)$. (The parallelotope will not be a
generating $r$-face unless it is the zonotope itself.)

A {\bf bounding face} is a translate of a generating face to the boundary of the zonotope using sums and differences
of generators not used in the definition of that face. For a generating facet
$\mathcal{F}=\mathcal{Z}(\boldsymbol{a}_{j_1},\ldots ,\boldsymbol{a}_{j_t})$, the associated bounding facets
can be given explicitly. Consider any $r-1$ linearly independent generators of the facet.
For example, suppose $\boldsymbol{a}_{j_1},\boldsymbol{a}_{j_2},\ldots, \boldsymbol{a}_{j_{r-1}}$ are independent.
The cross-product of these generators is then a normal vector to the facet. (See, for example, [4].)
We write this as $\boldsymbol{n}_{\mathcal{F}}=${\large $\times$}$({\mathitbf a}_{j_1},\ldots, {\mathitbf a}_{j_{r-1}})$.
Note that any two sets of $r-1$ independent generators of $\mathcal{F}$ will give cross-products that are scalar multiples
of each other. Relabel all generators $\boldsymbol{a}_1,\ldots ,\boldsymbol{a}_k$ of the zonotope as
$\boldsymbol{a}^0_1,\ldots,\boldsymbol{a}^0_{p},\boldsymbol{a}^{-}_{p+1},\ldots,
\boldsymbol{a}_q^{-},\,\boldsymbol{a}^{+}_{q+1},\ldots,\boldsymbol{a}^{+}_k$
with superscripts designating the generators with zero, negative, and
positive projections on $\boldsymbol{n}_\mathcal{F}$. It follows that
$\boldsymbol{a}^0_1,\ldots,\boldsymbol{a}^0_{p}$ is a maximal set of generators of rank $r-1$ with
$p=t$ and $\{\boldsymbol{a}^0_1,\ldots,\boldsymbol{a}^0_{p}\}=\{\boldsymbol{a}_{j_1},\ldots ,\boldsymbol{a}_{j_t}\}$, and that
$\mathcal{Z}(\boldsymbol{a}^0_1,\ldots,\boldsymbol{a}^0_{p})+\boldsymbol{a}^{-}_{p+1}+\cdots +\boldsymbol{a}^{-}_q$
will be one translation of $\mathcal{F}$ to a bounding facet, while
$\mathcal{Z}(\boldsymbol{a}^0_1,\ldots,\boldsymbol{a}^0_{p})+\boldsymbol{a}^{+}_{q+1}+\cdots +\boldsymbol{a}^{+}_k$ will be
the corresponding facet on the opposite side of the boundary.

\bigskip

We wish to revisit some results of Shephard and McMullen from [7-10] that examine how the central symmetry of the faces of a zonotope
relates to the symmetry of the entire zonotope. As zonotopes in their own right, the faces of a zonotope are always centrally symmetric.
For an arbitrary convex polytope, it turns out that the central symmetry of all faces of a given dimension implies the symmetry
of the faces of the next higher dimension, while the central symmetry of all faces in any dimension below that of the facets implies the
symmetry of the faces of the next lower dimension (McMullen, [7], [8]). Moreover, polytopes whose 2-faces are all centrally symmetric
are zonotopes. Consequently, zonotopes may be characterized as the convex polytopes of dimension $n$ whose faces of any one
particular dimension $k$ are centrally symmetric, where $2\leq k\leq n-2$.

In order to establish these and similar results, we start by considering zones of faces of polytopes.
Given a $k$-dimensional face $\mathcal{F}$ of polytope $\mathcal{P}$, the {\bf $k$-zone
$Z_k(\mathcal{F})$ induced by} $\mathcal{F}$ is defined as the union of all proper faces
that contain translates of $\mathcal{F}$ as faces. It clearly suffices to take the union only of facets, and
$Z_k(\mathcal{F})$ satisfies:

\medskip

\lcrline{}{\it if $k<j$, then
$Z_k(\mathcal{F})\,=\,\bigcup \{Z_j(\mathcal{F}^\prime)\,|\,\mathcal{F}\subset\mathcal{F}^\prime$ and $\mathcal{F}^\prime$ is a $j$-face$\}$.}{}

\medskip

\noindent The 1-zone $Z_1(\mathcal{E})\,=\,Z(\mathcal{E})$ induced by an edge $\mathcal{E}$ is called simply a {\bf zone}. It is the traditional
zone that give rise to the name zonotope.
\medskip

The following is a consequence of Shephard's Theorem 2, from [9].

\medskip

\noindent {\bf Lemma 1.2.}  {\it  Let $\mathcal{P}$ be a convex $n$-dimensional polytope in
$\mathbb{R}^n$ whose faces of dimension $(j+1)$ are all centrally symmetric,
where $(j+1)$ is such that $2\leq(j+1)\leq n$. Consider an orthogonal
projection of $\mathbb{R}^n$ to a complement of the $j$-dimensional affine
subspace supporting a particular $j$-dimensional face $\mathcal{F}$ of
$\mathcal{P}$. Then the image of $\mathcal{P}$ under this projection is an
$(n-j)$-dimensional convex polytope, $\pi(\mathcal{P})$, and all faces of
$\mathcal{P}$ of dimension $j$ that are translates or point reflection images of $\mathcal{F}$
map in one-to-one fashion to the vertices of $\pi(\mathcal{P})$. For each value $k$
with $j\leq k\leq n$, all $k$-dimensional faces of $\mathcal{P}$ containing
$\mathcal{F}$ are mapped in one-to-one fashion to all $($k$-$j$)$-dimensional
faces of $\pi(\mathcal{P})$ containing the image point of $\mathcal{F}$ under the projection.}

\medskip

Using this lemma, it is possible to give a new proof of an
$n$-dimensional version of a theorem of P. Alexandrov different from the proofs
given in [3] and [9].

\medskip

\noindent {\bf Proposition 1.3.}  {\it If all facets of a convex $n$-polytope ($n>2$)$\,$are centrally
symmetric, then the polytope is centrally symmetric.}

\smallskip

{\it Proof.} Consider an $n$-dimensional polytope $\mathcal{P}$ in $\mathbb{R}^n$.
Let $\mathcal{F}_1$ be a facet of $\mathcal{P}$ with center of symmetry
$\boldsymbol{c}_1$, and let $\mathcal{F}_{1,1}$ be an $(n-2)$-face of $\mathcal{P}$ that is
a facet of $\mathcal{F}_1$. Central symmetry ensures that the reflection
$(\mathcal{F}_{1,1})_{\boldsymbol{c}_1}=\mathcal{F}_{1,2}$ is the face of
$\mathcal{F}_1$ opposite to $\mathcal{F}_{1,1}$. This face is shared with an
adjacent facet, $\mathcal{F}_2$. Let $\boldsymbol{c}_2$ be the center of
$\mathcal{F}_2$. The reflection
$(\mathcal{F}_{1,2})_{\boldsymbol{c}_2}=\mathcal{F}_{1,3}$ is then the face
opposite $\mathcal{F}_{1,2}$ on the boundary of $\mathcal{F}_2$. (It is also
a translate of $\mathcal{F}_{1,1}$.) Face $\mathcal{F}_{1,3}$ is shared
with another facet, $\mathcal{F}_3$. In this way, successive $(n-2)$-faces
$\mathcal{F}_{1,1},\mathcal{F}_{1,2},\ldots,\mathcal{F}_{1,m_1+1}=%
\mathcal{F}_{1,1}$ are determined that are alternately point reflections and
translations of $\mathcal{F}_{1,1}$. The faces determine a corresponding
chain of facets,
$\mathcal{F}_1,\mathcal{F}_2,\ldots,\mathcal{F}_{m_1+1}=\mathcal{F}_1$, whose
union, $\mathcal{F}_1\cup\mathcal{F}_2\cup\cdots\cup\mathcal{F}_{m_1}=Z_{(n-2)}(\mathcal{F}_{1,1})$,
is an $(n-2)$-zone on the boundary of $\mathcal{P}$.

Choose an $(n-2)$-face $\mathcal{F}_{2,1}$ adjacent to $\mathcal{F}_{1,1}$ on the
boundary of $\mathcal{F}_1$. This determines another sequence of
$(n-2)$-dimensional faces,
$\mathcal{F}_{2,1},\mathcal{F}_{2,2},\ldots,\mathcal{F}_{2,m_2+1}=%
\mathcal{F}_{2,1}$, consisting of reflected and translated copies of
$\mathcal{F}_{2,1}$ and another sequence of facets
$\mathcal{F}'\!\!_1,\mathcal{F}'\!\!_2,\ldots,\mathcal{F}'\!\!_{m_2+1}=\mathcal{F}'\!\!_1$
whose union is a second $(n-2)$-zone, $Z_{(n-2)}(\mathcal{F}_{2,1})$, on the boundary of $\mathcal{P}$.
Project $\mathbb{R}^n$ to the orthogonal complement of the $(n-3)$-dimensional affine
subspace supporting the face $\mathcal{F}_{1,2,1}=\mathcal{F}_{1,1}\cap\mathcal{F}_{2,1}$.
It follows from Lemma 1.2 that the projections of the two $(n-2)$-zones of facets become zones of
$2$-faces on the boundary of 3-dimensional $\pi(\mathcal{P})$. Zones on a convex polyhedron
are circumferential; any two intersect twice. As the projected zones on
$\pi(\mathcal{P})$ both include $\pi(\mathcal{F}_1)$, they must
therefore intersect a second time. It follows that the $(n-2)$-zones of preimages
must also intersect twice. In other words, if $(n-2)$-zones on the boundary of a
convex $n$-dimensional polytope with centrally symmetric facets intersect at
all, then they intersect twice.

The two $(n-2)$-zones under consideration intersect at
$\mathcal{F}_1^{\,\prime}=\mathcal{F}_1$. Hence they also intersect at
$\mathcal{F}^{\,\prime}_j=\mathcal{F}_k$ for some $j,k>1$. Facet
$\mathcal{F}_k$ then includes translative or reflective copies of both $\mathcal{F}_{1,1}$
and $\mathcal{F}_{2,1}$ as part of its boundary. The same is true for
$\mathcal{F}_1$. The two facets are therefore parallel. By convexity,
$\mathcal{F}_k$ is the unique facet of $\mathcal{P}$ parallel to
$\mathcal{F}_1$. Denote it as $\mathcal{F}_1^{\text{ }\,\text{op}}$. In this
way, every facet of $\mathcal{P}$ is paired with a unique parallel, opposite
facet. In particular, $\mathcal{P}$ and all zones of facets of $\mathcal{P}$
contain even numbers of facets in parallel, opposite pairs.

We wish to show that $\mathcal{F}_1$ and $\mathcal{F}_1^{\,\text{op}}$ are point symmetric
images of each other. To see that this is so, let $\boldsymbol{c}_1$ be the
center of $\mathcal{F}_1$ and let $\boldsymbol{c}_1^{\text{op}}$ be the
center of $\mathcal{F}_1^{\,\text{op}}$. By assumption, $\boldsymbol{c}_1^{\text{op}}$ exists.
Set $\boldsymbol{c}=\frac{1}{2}(\boldsymbol{c}_1+\boldsymbol{c}_1^{\,%
\text{op}})$. Consider the symmetric image $(\mathcal{F}_1)_{\boldsymbol{c}}$
of $\mathcal{F}_1$, which is a centrally symmetric $(n-1)$-polytope that must
lie in the same hyperplane, $\mathcal{H}$, as $\mathcal{F}_1^{\,\text{op}}$.
We will see that $(\mathcal{F}_1)_{\boldsymbol{c}}$ and
$\mathcal{F}_1^{\,\text{op}}$ are identical. Each can be defined in terms of
the intersection of $\mathcal{H}$ with half-spaces determined by the
hyperplanes supporting all adjacent facets. Suppose $\mathcal{F}_\ast$ is a
facet of $\mathcal{P}$ adjacent to $\mathcal{F}_1$ with center $\boldsymbol{c}_\ast$. Then
$(\mathcal{F}_\ast)_{\boldsymbol{c}}$ will be adjacent to
$(\mathcal{F}_1)_{\boldsymbol{c}}\,$, and $\mathcal{F}_\ast^{\,\text{op}}$ will
be adjacent to $\mathcal{F}_1^{\,\text{op}}$. Denote the hyperplanes supporting
$\mathcal{F}_\ast$, $(\mathcal{F}_\ast)_{\boldsymbol{c}}$, and
$\mathcal{F}_\ast^{\,\text{op}}$ by $\mathcal{H}_\ast,$
$(\mathcal{H}_\ast)_{\boldsymbol{c}},$ and $(\mathcal{H}_\ast)^{\text{op}}$
respectively. These hyperplanes are parallel, and the latter two contain the
center of symmetry
$(\boldsymbol{c}_{_{\scriptstyle\ast}})_{\boldsymbol{c}}=\boldsymbol{c}_{%
\scriptstyle\ast}^{\text{op}}$ common to both
$(\mathcal{F}_\ast)_{\boldsymbol{c}}$, and $\mathcal{F}_\ast^{\,\text{op}}$.
Hence $(\mathcal{H}_\ast)_{\boldsymbol{c}}=(\mathcal{H}_\ast)^{\text{op}}$.
This hyperplane defines two half-spaces one of which contains $\mathcal{F}_\ast$ and
is included among the half-spaces whose intersection with $\mathcal{H}$
defines both $(\mathcal{F}_1)_{\boldsymbol{c}}$, and
$\mathcal{F}_1^{\,\text{op}}$. The other half-spaces defining the two
facets are determined in a similar manner. As a result,
$(\mathcal{F}_1)_{\boldsymbol{c}}$, and $\mathcal{F}_1^{\,\text{op}}$ have
the same definition in terms of intersections, and so
$(\mathcal{F}_1)_{\boldsymbol{c}}=\mathcal{F}_1^{\,\text{op}}$. Consequently,
\ $\mathcal{F}_1$ and $\mathcal{F}_1^{\,\text{op}}$ are point reflections of each
other. $\mathcal{F}_1\cup\mathcal{F}_1^{\,\text{op}}$ is therefore centrally
symmetric with center of symmetry
$\boldsymbol{c}_{1,1}\,{\mathrel{\mathop :}=}\,\,\boldsymbol{c}\,=\,\frac{1}{2}(\boldsymbol{c}_1+%
\boldsymbol{c}_1^{\,\text{op}})$. The same can then be said for all facets of
$\mathcal{P}$.

Thus, for each facet $\mathcal{F}_j$ of $\mathcal{P}$, the union
$\mathcal{F}_j\cup\mathcal{F}_j^{\,\text{op}}$ is centrally symmetric with
center of symmetry
$\boldsymbol{c}_{j,j}\,\,{\mathrel{\mathop :}=}\,\,\frac{1}{2}(\boldsymbol{c}_j+\boldsymbol{c}_j^{\,%
\text{op}})$. Moreover, if $\mathcal{F}_j$ and $\mathcal{F}_k$ are adjacent
facets sharing an $(n-2)$-face, the centers of symmetry agree on that face
making those centers the same: $\boldsymbol{c}_{j,j}=\boldsymbol{c}_{k,k}$.
By moving around the entire boundary of $\mathcal{P}$ from facet to adjacent
facet, all opposite pairs of facets share a common center of symmetry. This
becomes the center of symmetry for the entire polytope, which is therefore
centrally symmetric. \hfill{\qedsymbol}

\medskip

An easy inductive argument then gives the following immediate consequence:

\bigskip

\noindent {\bf Corollary 1.4.} {\it If the $k$-dimensional faces of
an $m$-dimensional convex polytope in $\mathbb{R}^n$ are centrally
symmetric for some value $k\geq 2$, then the $(k+1)$-dimensional
faces are also centrally symmetric.}

\medskip

\noindent Proposition 1.3 and Corollary 1.4 apply when $k\geq 2$. When
$k=1$ and $m=2$, a convex polygon has 1-dimensional edges
that are centrally symmetric, but the polygon itself need not be centrally
symmetric. The following conditions show when an arbitrary polygon, or when
any closed configuration of line segments, is centrally symmetric.

\medskip

\noindent {\bf Proposition 1.5.} {\it A $2$-dimensional convex polygon is centrally
symmetric if and only if it has an even number of edges and all pairs of opposite edges
are parallel and of equal length. More generally, a closed configuration consisting of an even
number of directed line segments $\boldsymbol{s}_1,\ldots,\boldsymbol{s}_{2t}\subset\mathbb{R}^n$
with the end point of each segment coinciding with the starting point of the next considered modulo $2t$
is centrally symmetric if and only if for each $j=1,\ldots,t,\,\boldsymbol{s}_j$ and
$\boldsymbol{s}_{t+j}$ are parallel, of equal length, and of opposite orientation.}

\smallskip

{\it Proof.} For a convex polygon with an even number of edges, if a direction chosen for one edge is
used to determine a consistent direction for all successive edges, and if opposite edges are always equal
and parallel, then all conditions for a closed configuration of directed line segments given in the
statement of the proposition will be satisfied.

Suppose that for such a configuration, $\boldsymbol{s}_1,\ldots,\boldsymbol{s}_{2t}$, each pair of
opposite segments---those of the form $\boldsymbol{s}_j$ and $\boldsymbol{s}_{t+j}$---are parallel,
equal, and of opposite orientation. The segments from a pair then determine a parallelogram in the plane
they span. Their opposite orientation ensures that one diagonal of the parallelogram will
connect the starting points of the segments, and the other diagonal will connect the ending points.
The diagonals cross and are divided in halves at the center of the parallelogram. This center and either one
of the pair of oppositely oriented segments define a symmetric cone that connects the segments as symmetric
images of each other. The diagonals of the parallelogram serve as the extreme, bounding elements of the cone.
Symmetric cones of adjacent pairs of segments from the configuration will share one or the
other of these bounding diagonal elements, and hence will have the same centers of symmetry.
In this way, the centers of symmetry of all opposite pairs of directed segments
in the configuration will be the same, and hence the configuration will be centrally symmetric.

Conversely, if the configuration is centrally symmetric, then reflection of any $\boldsymbol{s}_j$
through the center of symmetry will produce the opposite segment
$\boldsymbol{s}_{t+j}$, which must then be parallel to $\boldsymbol{s}_j$, of equal length,
and with opposite orientation. Consequently, these conditions are equivalent to central symmetry
of the configuration. \hfill{\qedsymbol}

\medskip

\noindent {\bf Proposition 1.6.} {\it A convex $m$-polytope in $\mathbb{R}^n$ in which all
$2$-faces are centrally symmetric is a zonotope.}

\smallskip

{\it Proof.} When $m=2$, consider a centrally symmetric filled-in polygon in the
plane. Select any edge of the polygon together with the opposite edge, which
is its reflection with respect to the center of symmetry. Connecting the
endpoints of the two edges produces a strip in the form of a parallelogram,
which by virtue of convexity is wholly contained within the polygon. If the
strip is removed, translation by a vector defined by the selected edge joins
the endpoints of each removed edge to produce a smaller, centrally symmetric
polygon with two fewer edges. The original polygon is the Minkowski sum of
the new polygon and the selected edge. By downward induction on the number of edges,
the original polygon becomes the Minkowski sum of edges and hence is a zonogon.

Now suppose $m\geq 3$. Assume as an induction hypothesis that the proposition
is true for all polytopes of dimension$\,\,<m$. Consider an $m$-dimensional
convex polytope $\mathcal{P}$ with centrally symmetric 2-faces. Select an
edge $\mathcal{E}$ of $\mathcal{P}$ and let the zone
$Z(\mathcal{E})=\mathcal{F}_1\cup\mathcal{F}_2\cup\cdots\cup\mathcal{F}_t$
be the union of all facets containing translates of $\mathcal{E}$ as edges. All facets
belonging to $Z(\mathcal{E})$ have centrally symmetric 2-faces, so by
the induction hypothesis, each
is a zonotope containing $\mathcal{E}$ as an edge.
Each $\mathcal{F}_j$ therefore decomposes as the Minkowski sum of $\mathcal{E}$ with a
smaller zonotope, $\mathcal{F}^\ast_j$. Thus,
$\mathcal{F}_j=\mathcal{F}^\ast_j\,{\oplus}\,\mathcal{E}$. It follows that

\medskip

\lcrline{}{$Z(\mathcal{E})\,=\,\mathcal{F}_1\cup\mathcal{F}_2\cup\cdots\cup\mathcal{F}_t=(%
\mathcal{F}^\ast_1\,{\oplus}\,\mathcal{E})\cup(\mathcal{F}^%
\ast_2\,{\oplus}\,\mathcal{E})\cup\cdots\cup(\mathcal{F}^%
\ast_t\,{\oplus}\,\mathcal{E})=(\mathcal{F}^\ast_1\cup%
\mathcal{F}^\ast_2\cup\cdots\cup\mathcal{F}^\ast_t)\,{\oplus}\,\mathcal{E}$.}{}

\medskip

\noindent(Note that for any subsets $A,B,C\subset\mathbb{R}^n$,
$(A\,{\oplus}\,C)\cup(B\,{\oplus}\,C)=(A\cup
B)\,{\oplus}\,C.)\,$\,\,Suppose the remaining facets of
$\mathcal{P}$ are
$\mathcal{F}^{\,\prime}_1,\mathcal{F}^{\,\prime}_2,\ldots,\mathcal{F}^{\,%
\prime}_s$ so that the boundary is
$\,(\mathcal{F}_1\cup\mathcal{F}_2\cup\cdots\cup\mathcal{F}_t)\cup(%
\mathcal{F}^{\,\prime}_1\cup\mathcal{F}^{\,\prime}_2\cup\cdots\cup%
\mathcal{F}^{\,\prime}_s)$. After replacing each $\mathcal{F}_j$ with
$\mathcal{F}^\ast_j$, it follows that

\smallskip

\lcrline{}{$(\mathcal{F}^\ast_1\cup\mathcal{F}^\ast_2\cup\cdots\cup%
\mathcal{F}^\ast_t)\cup(\mathcal{F}^{\,\prime}_1\cup\mathcal{F}^{\,\prime}_2%
\cup\cdots\cup\mathcal{F}^{\,\prime}_s)$}{}

\medskip

\noindent is the boundary of a polytope $\mathcal{P}^{\prime}$ with fewer edges than
$\mathcal{P}$, with centrally symmetric 2-faces, and with
$\mathcal{P}=\mathcal{P}^{\prime}\,{\oplus}\,\mathcal{E}$.
Once more, by downward induction on the number of non-parallel edges,
$\mathcal{P}$ will become the Minkowski sum of edges and hence a zonotope. In
this way, every convex polytope of dimension $m$ with centrally symmetric
2-faces will be a zonotope.\hfill{\qedsymbol}

\bigskip

\noindent A different proof of Proposition 1.6 can be found in  [3, Proposition 2.2.14].

\medskip

\noindent {\bf Corollary 1.7.} {\it A convex polytope is a zonotope if and only if it decomposes
into zonotopes. Equivalently, a convex polytope is a zonotope if and only if
it decomposes into parallelotopes.}

\smallskip

{\it Proof.} A zonotope trivially decomposes into zonotopes. By a theorem of
Shephard and McMullen (see [4] or [10]),
it decomposes into parallelotopes. (Note that the decomposition, also called a {\bf tiling}, means
that the zonotope is the union of parallelotopes meeting each other in lower-dimensional facets.)

Conversely, suppose an $m$-polytope $\mathcal{P}\subset\mathbb{R}^n$
decomposes into zonotopes, and hence into parallelotopes. It follows that in
every dimension$\,\,<m$, every face of $\mathcal{P}$ also decomposes into
parallelotopes. In particular, each 2-face decomposes into (filled-in)
parallelograms and each edge decomposes into edges from those parallelograms.
Consider a specific 2-face $\mathcal{F}$ of $\mathcal{P}$, a particular edge of $\mathcal{F}$,
and a parallelogram $\mathcal{P}^2$ that is
part of the decomposition of $\mathcal{F}$ and has an edge contained in the
designated edge of $\mathcal{F}$. The edge of $\mathcal{P}^2$ opposite to
the one lying in the edge of $\mathcal{F}$ is itself shared with another
parallelogram in the decomposition of $\mathcal{F}$. By tracking this edge
from parallelogram to parallelogram, a strip of parallelograms sharing
translated copies of the edge extends across $\mathcal{F}$ to its
far side. The process can be repeated with each of the parallelograms that
shares an edge with part of the designated edge of $\mathcal{F}$. Taken
together, the resulting strips produce a translated copy of the designated
edge of $\mathcal{F}$ on the far side of the boundary of $\mathcal{F}$. (A
complete edge of $\mathcal{F}$ must be obtained in this way because if part
of an edge on the far side was not reached by such a strip of parallelograms,
a strip formed in reverse would produce a copy of that part back on the
originally designated edge of $\mathcal{F}$.) In this way, every edge of
$\mathcal{F}$ is paired with a translated, parallel, opposite copy of that
edge. It then follows from Proposition 1.5 that $\mathcal{F}$ is centrally
symmetric.

Once all 2-faces are centrally symmetric, the polytope is a zonotope by
Proposition 1.6.

\hfill{\qedsymbol}

\smallskip

McMullen [7] demonstrated that central symmetry for faces
migrates to lower as well as higher dimensions in a convex polytope
provided one starts by assuming the central symmetry of faces in a dimension
lower than that of the facets. We give McMullen's proof
rephrased in the current notation.

\medskip

\noindent {\bf Proposition 1.8.} {\it If the $(n-2)$-dimensional faces of $n$-dimensional convex
polytope $\mathcal{P}\subset\mathbb{R}^n$ are centrally symmetric, then the
$(n-3)$-dimensional faces are centrally symmetric.}

\smallskip

{\it Proof.}  Consider an $(n-3)$-face, $\mathcal{F}_{1,1,1}$, on the boundary of
$(n-2)$-face, $\mathcal{F}_{1,1}$, which is in turn on the boundary of facet,
$\mathcal{F}_1$, of $\mathcal{P}$. Central symmetry implies there are $(n-3)$-zones
of $(n-2)$-faces on the boundary of $\mathcal{F}_1$ induced by $\mathcal{F}_{1,1,1}$.

If an $(n-3)$-zone of $(n-2)$-faces induced by some $(n-3)$-face is of length four, then
the $(n-3)$-face must be centrally symmetric. To see this,
suppose $\mathcal{F}_{1,1}\cup\mathcal{F}_{1,2}\cup\mathcal{F}_{1,3}\cup%
\mathcal{F}_{1,4}$ is a zone of length four induced by $\mathcal{F}_{1,1,1}$. From Proposition 1.3,
the facet $\mathcal{F}_1$ is centrally symmetric, so the face
opposite $\mathcal{F}_{1,1}$ in this zone satisfies
$\mathcal{F}_{1,3}=(\mathcal{F}_{1,1\scriptstyle})_{\boldsymbol{c}_1}$ where $\boldsymbol{c}_1$ is the
center of $\mathcal{F}_1$. In particular,

\lcrline{}{$\mathcal{F}_{1,2}\cap\mathcal{F}_{1,3}=(\mathcal{F}_{1,1,1})_{%
\scriptstyle\boldsymbol{c}_1}$.}{}

\smallskip

\noindent At the same time, central symmetry of $\mathcal{F}_{1,1}$ followed by central
symmetry of $\mathcal{F}_{1,2\scriptstyle}$ imply by Lemma 1.1$(a)$ that

\lcrline{}{$\mathcal{F}_{1,2}\cap\mathcal{F}_{1,3}=2(\boldsymbol{c}_{1,2}-%
\boldsymbol{c}_{1,1})+\mathcal{F}_{1,1,1}$.}{}

\medskip

\noindent From these two equations, Lemma 1.1$(e)$ implies that $\mathcal{F}_{1,1,1}$ is
centrally symmetric.

To complete the proof, it suffices to demonstrate that any $(n-3)$-face such
as $\mathcal{F}_{1,1,1}$ must be contained in some $(n-3)$-zone of
length four on the boundary of $\mathcal{F}_{1}$. This can be done by projecting
$\mathbb{R}^n$ orthogonally to the complement of the affine $(n-3)$-dimensional
subspace containing $\mathcal{F}_{1,1,1}$. The image of $\mathcal{P}$ under this
projection is a 3-dimensional centrally symmetric polyhedron with centrally symmetric
facets---a zonohedron, $\pi(\mathcal{P})$. Lemma 1.2 guarantees that images of
translated and reflected copies of $\mathcal{F}_{1,1,1}$ remain distinct on the
boundary of $\pi(\mathcal{P})$ and become its vertices, and that the images of the
$(n-2)$-faces containing $\mathcal{F}_{1,1,1}$ become, in one-to-one fashion,
the edges on the boundary of $\pi(\mathcal{P})$. The existence of a
parallelogram of edges bounding a face of $\pi(\mathcal{P})$ would therefore
demonstrate that the preimage is an $(n-3)$-zone of $(n-2)$-faces containing
$\mathcal{F}_{1,1,1}$ of length four on the boundary of $\mathcal{F}_1$. But
such a parallelogram must exist on the boundary of $\pi(\mathcal{P})$ because
any zonohedron contains at least six parallelogram faces, as can be seen
using Euler's formula.

In more detail, if $f_j=\,\,$\#$\,$faces with $j$ edges and $v_j=\,\,$\#$\,$vertices
of valency $j$, then $f=\sum f_j,\,v=\sum v_j,\!$ and $2e=\sum jf_j=\sum jv_j$.
Using these values in Euler's formula$\,$ $v-e+f=2\,$ yields two equations,

\medskip

\lcrline{$\displaystyle$}{$\sum(2-j)f_j+\sum 2v_j=4$ \ and \ $\sum
2f_j+\sum(2-j)v_j=4$.}{}

\smallskip

\noindent Combining the first equation with twice the second,

\medskip

\lcrline{$\displaystyle$}{$\sum_{j\geq 3}(6-j)f_j+\sum_{j\geq
3}(6-2j)v_j=12,$}{}

\smallskip

\noindent from which it follows that

\smallskip

\lcrline{$\displaystyle$}{$3f_3+2f_4+f_5\geq 12+\sum_{j\geq 7}(j-6)f_j$.}{}

\medskip

\noindent When all faces are centrally symmetric, $f_3=f_5=0$, so from the preceding inequality, $f_4\geq 6$.

\hfill{\qedsymbol}

\bigskip

\bigskip

Propositions 1.3 and 1.8 can be combined to obtain:

\bigskip

\noindent {\bf Theorem 1.9.} {\it For a polytope of dimension $m$ in $\mathbb{R}^n$, $m\leq n$,
if all $j$-dimensional faces are centrally symmetric for a particular value,
$2 \leq j \leq(m-2)$, then the faces in every dimension, including the polytope
itself, are centrally symmetric and the polytope is a zonotope.}

\bigskip

A point made in the proof of Proposition 1.8 is that existence of a $k$-zone of length
$\equiv 0\,$mod$\,4$ implies that the $k$-face generating the zone is centrally symmetric.
If one considers zones of a specific polytope, the possible lengths for the zones
are limited by the nature of the $k$-faces. For example,
the 24-cell $\mathcal{P}_{24}\subset\mathbb{R}^4$ is a regular polytope
with twenty four facets, each of which is a regular
octahedron. The four pairs of opposite $2$-faces of a particular
facet give rise to four $2$-zones of length $6\equiv 2\,$mod$\,4$
on the boundary. No zone can have a length that is a multiple of 4 because
the generating 2-face of such a zone would have to be centrally symmetric, which is not
the case for the triangular $2$-faces of this polytope.

The $2$-zones of $\mathcal{P}_{24}$ also provide a discrete Hopf fibration of its boundary.
Starting from a particular facet $\mathcal{F}_1$, the four zones of facets generated
by opposite pairs of $2$-faces of this facet
are each of length 6. The zones meet at $\mathcal{F}_1$ and again in their
fourth facets, denoted $\mathcal{F}_1^{\text{op}}$. The zones can be written as
$\mathcal{F}_1\cup\mathcal{F}^j_2\cup\mathcal{F}_3^j\cup\mathcal{F}_1^{%
\text{op}}\cup\mathcal{F}_5^j\cup\mathcal{F}_6^j$ for $j=1,\ldots,4$.
Together, these zones account for $(4\cdot 4)+2=18$ of the facets. The remaining six
facets, labeled $\mathcal{F}_1^\ast,\ldots,\mathcal{F}_6^\ast$, fill the interstices
and complete the boundary of $\mathcal{P}_{24}$. Consider one of the zones, say
$\mathcal{F}_1\cup\mathcal{F}^1_2\cup\mathcal{F}_3^1\cup\mathcal{F}_1^{%
\text{op}}\cup\mathcal{F}_5^1\cup\mathcal{F}_6^1$. This zone and three new $2$-zones
$\mathcal{F}_2^2\cup\mathcal{F}_1^\ast\cup\mathcal{F}_3^3\cup\mathcal{F}_5^4%
\cup\mathcal{F}_2^\ast\cup\mathcal{F}_6^3$,
$\mathcal{F}_2^3\cup\mathcal{F}_3^\ast\cup\mathcal{F}_3^4\cup\mathcal{F}_5^2%
\cup\mathcal{F}_4^\ast\cup\mathcal{F}_6^4$, and
$\mathcal{F}_2^4\cup\mathcal{F}_5^\ast\cup\mathcal{F}_3^2\cup\mathcal{F}_5^3%
\cup\mathcal{F}_6^\ast\cup\mathcal{F}_6^2$
are mutually disjoint. Together, they include all the facets of $\mathcal{P}_{24}$ and
constitute a discrete Hopf fibration of the boundary of $\mathcal{P}_{24}$.
Other examples with similar fibrations
are two more regular polytopes in $\mathbb{R}^4$: the 120-cell,
which has regular dodecahedral facets, and the 600-cell, which has tetrahedral facets.
Prisms in $\mathbb{R}^4$ also have simple discrete fibrations. This raises the question of whether
there might exist sequences of polytopes in $\mathbb{R}^4$ with increasing numbers of $2$-zones that
allow discrete Hopf fibrations, which in the limit give the fibration of the 3-sphere.

\medskip

With regard to $(n-2)$-zones---as opposed to $1$-zones---on the boundaries of zonotopes, start by considering
the $4$-cube. One standard 3-dimensional projection of the $4$-cube consists of
inner and outer cubes whose corresponding vertices are connected by additional edges. The facets in this
projection consist of two 3-cubes that can be labeled {\it inner} and {\it outer}, and six more $3$-cells
surrounding the inner cube that can be labeled in pairs as {\it up/down}, {\it front/back},
and {\it left/right}. Two non-intersecting $(n-2)=2$-zones consisting of the facets
{\it front-down-back-up} and {\it inner-left-outer-right} then form a decomposition
of the boundary of the 4-cube. This is the simplest discrete version of the Hopf fibration
$S^1\hookrightarrow S^3\rightarrow S^2$ for spheres that is realizable for a polytope.
\medskip

While some $(n-2)$-zones on the boundaries of $n$-zonotopes
might not intersect, an argument given in the course of the proof of Proposition 1.3 establishes

\medskip

\noindent {\bf Proposition 1.11.} {\it If two $(n-2)$-zones on the boundary of a zonotope intersect, then they intersect
precisely twice.}

\section*{2. \,\,\,\,Congruences of Zonotopes}

For the study of congruence, start with an identity that comes from the matrix
$A^TA$, which will be called the {\bf shape matrix} of ${\mathcal{Z}(A)}$. If two zonotopes have
the same shape matrix, they are congruent because the transformation taking
generating vectors of one zonotope to corresponding
vectors of the other is an isometry. The matrix formulation is the following:

\bigskip

\noindent {\bf Proposition 2.1.}  {\it If $A$ and $B$ are $n\times k$ matrices ($n$ and $k$
arbitrary), then $A^TA=B^TB$ if and only if $B=QA$ where $Q$ is an
$n\times n$ orthogonal matrix.}

\medskip

{\it Proof.} We prove only the non-trivial direction and assume $A^TA=B^TB$. Let
$A=[\boldsymbol{a}_1,\ldots,\boldsymbol{a}_k]$ and
$B=[\boldsymbol{b}_1,\ldots,\boldsymbol{b}_k]$.  Observe first that
independence of the columns of $A$ is equivalent to nonsingularity of $A^TA$.
(Independence of the columns and
$A^TA\boldsymbol{x}=\boldsymbol{0}\Rightarrow\boldsymbol{x}^TA^TA%
\boldsymbol{x}=\boldsymbol{0}\Rightarrow\left\|A\boldsymbol{x}\right\|^2=0%
\Rightarrow\,A\boldsymbol{x}=\boldsymbol{0}\,\!\Rightarrow\boldsymbol{x}=%
\boldsymbol{0}$.  Nonsingularity of $A^TA$ and
$A\boldsymbol{x}=\boldsymbol{0}\Rightarrow
A^TA\boldsymbol{x}=\boldsymbol{0}\,\Rightarrow\boldsymbol{x}=%
\boldsymbol{0}$.) The same can be said for $B$.

{\bf Case 1:} $A$ has independent columns. By the initial observation,
we have independence of the columns of $A$ iff $A^TA=B^TB$ is
nonsingular, that is, iff we have independence of the columns of $B$. Hence
col$(A)^\perp$ and col$(B)^\perp$ both have$\,$dimension $n-k$. Let
$\boldsymbol{a}_{k+1},\ldots,\boldsymbol{a}_n$ and
$\boldsymbol{b}_{k+1},\ldots,\boldsymbol{b}_n$ be orthonormal bases of
col$(A)^\perp$ and col$(B)^\perp$ respectively. Let $Q$ be the matrix of
the transformation defined by $Q\boldsymbol{a}_i=\boldsymbol{b}_i$ for
$i=1,\ldots,n$. Thus, in particular, $QA=B$. Clearly, $Q$ is orthogonal
on col$(A)^\perp$. It then remains to be shown that $Q$ is also orthogonal on
col$(A)$. To do so, consider vectors
$\boldsymbol{x},\boldsymbol{y}\in$col$(A)$ written as
$\boldsymbol{x}=A\boldsymbol{c}$ and $\boldsymbol{y}$ $=A\boldsymbol{d}$.
We then have

\smallskip

\lcrline{}{$(Q\boldsymbol{x})\cdot(Q\boldsymbol{y})=(QA\boldsymbol{c})^T(QA%
\boldsymbol{d})$ $\!=(B\boldsymbol{c})^T(B\boldsymbol{d})=$
$\boldsymbol{c}^TB^TB\boldsymbol{d}=\boldsymbol{c}^TA^TA\boldsymbol{d}=(A%
\boldsymbol{c})^T(A\boldsymbol{d})=\boldsymbol{x}\cdot\boldsymbol{y}$.}{}

\smallskip

\noindent Hence $Q$ preserves inner products on col$(A)$ and is therefore orthogonal.

\smallskip

{\bf Case 2:} $A$ has dependent columns. Let
$A_0=\left[\boldsymbol{a}_{i_1},\ldots,\boldsymbol{a}_{i_l}\right]$ $(l<k)$
where $\boldsymbol{a}_{i_1},\ldots,\boldsymbol{a}_{i_l}$ is a maximal
collection of independent columns of $A$, and set
$B_0=\left[\boldsymbol{b}_{i_1},\ldots,\boldsymbol{b}_{i_l}\right]$. The
hypothesis $A^TA=B^TB$ implies $A_0^TA_0=B_0^TB_0$, so from case 1,
$B_0=QA_0$ for some orthogonal $Q$. We wish to show $B=QA$ and so that
$Q\boldsymbol{a}_t=\boldsymbol{b}_t$ for each $t\neq i_1,\ldots,i_l$. By
the observation made at the start of the proof, maximal independence of the
columns of $A_0$ implies the same for the columns of $B_0$. As
$\boldsymbol{a}_t$ and $\boldsymbol{b}_t$ are thus dependent on the columns
of $A_0$ and $B_0$ respectively, we may write
$\boldsymbol{a}_t=\,{\sum_{j=1}^lc_j\boldsymbol{a}_{i_j}}=A_0\boldsymbol{c}$
and $\boldsymbol{b}_t=\sum_{j=1}^ld_j\boldsymbol{b}_{i_j}=B_0\boldsymbol{d}$.
\ Meanwhile, $Q\boldsymbol{a}_t=QA_0\boldsymbol{c}=B_0\boldsymbol{c}$. To
complete the proof, we must show $\boldsymbol{c}=\boldsymbol{d}$. But this
follows from the series of implications:

\smallskip

\lcrline{}{\small{$A^TA=B^TB\Rightarrow
A_0^T\boldsymbol{a}_t=B_0^T\boldsymbol{b}_t\Rightarrow
A_0^TA_0\boldsymbol{c}=B_0^TB_0\boldsymbol{d}=A_0^TA_0\boldsymbol{d}%
\Rightarrow
A_0^TA_0(\boldsymbol{c}-\boldsymbol{d})=\boldsymbol{0}\Rightarrow%
\boldsymbol{c}-\boldsymbol{d}=\boldsymbol{0}$.}}{}

\smallskip

\noindent The last implication follows from the nonsingularity of $A_0^TA_0$, which
itself follows from yet another application of the initial observation of the
proof.\hfill{\qedsymbol}

\bigskip

Applications of this proposition (in the complex case and with a different
proof) were given in [6], but no geometric interpretations
involving zonotopes were mentioned. Some results
from that article take on added significance when the matrices (in the
real case) are regarded as shape matrices of zonotopes. The proposition in
the form given here together with some of its consequences represent past
joint work with Nishan Krikorian. We now extend some of the results that
relate to zonotopes.

One piece of information that the shape matrix does not contain is
the dimension of the space in which the zonotope resides.
What if two such objects have the same shape matrix but lie in different
dimensional Euclidean spaces? Then Proposition 2.1 becomes:

\medskip

{\it If $A$ is $m\times k$ and $B$ is $n\times k$ ($m\leq n$), then
$A^TA=B^TB$ iff $B=QA$ where $Q\,is$ $n\times m$ with orthonormal
columns.}

\medskip

\noindent This is just a slight generalization whose proof is omitted.

\medskip

Another piece of geometric information about a zonotope that the shape matrix
does not give is its embedding in Euclidean space.
So far, congruent zonotopes have implicitly been assumed attached to the origin
at vertices that correspond to each other under the congruence. But if that is
not the case, such a correspondence can be made after altering the defining matrix
of one of the zonotopes in order to change the vertex located at the origin.
For example, if $\mathcal{Z}$($\boldsymbol{a}_1,\ldots,\boldsymbol{a}_k$) is
translated along edge $\boldsymbol{a}_1$ so that the origin is moved to the
terminal point of that edge, the resulting copy of the zonotope will have the form
$-\boldsymbol{a}_1+\mathcal{Z}$($\boldsymbol{a}_1,\ldots,\boldsymbol{a}_k$)$\,%
=\mathcal{Z}$($-\boldsymbol{a}_1,\boldsymbol{a}_2,\ldots,%
\boldsymbol{a}_k$). A sequence of such translations can be used to reach
any vertex yielding $\mathcal{Z}$($\boldsymbol{a}'_1,\ldots,\boldsymbol{a}'_k$),
which will therefore be related to
$\mathcal{Z}$($\boldsymbol{a}_1,\ldots,\boldsymbol{a}_k$) simply by
$\boldsymbol{a}'_i=\pm\boldsymbol{a}_i$. The generating matrices will then be
related by $A'=AJ$ and the shape matrices by $(A')^TA'=JA^TAJ$ for
some $k\times k$ diagonal matrix $J$ with $\pm 1$'s on the principal
diagonal. If, in addition, we wish to reorder the generating vectors (to
match, for example, the order of generating vectors of some congruent
zonotope), this can be done by pre-multiplying $A$ by a
permutation matrix, $\varSigma$. The general statement about congruence is
then:

\medskip

\noindent {\bf Theorem 2.2.} {\it $\mathcal{Z}(A)$ and $\mathcal{Z}(B)$ are congruent,
where $A$ is $m\times k$ and $B$ is $n\times k$ ($m\leq n$),
\text{if and only if} $($A$')^TA'=B^TB$ where $A'=A\varSigma J$ for some
$k\times k$ permutation matrix $\varSigma$ and some diagonal matrix $J$
with $\pm 1$'$s\,on\,the\,diagonal$, or equivalently, if and only if there
exists an $n\times m$ matrix $Q$ with orthonormal columns such that
$B=QA\varSigma J$}.

\medskip

Now, consider a pair of generating matrices $A$ and $B$ of size $n\times k$
with independent columns, along with the parallelotopes $\mathcal{P}(A)$ and
$\mathcal{P}(B)$ in $\mathbb{R}^n$, and the zonotopes $\mathcal{Z}(A^T)$ and
$\mathcal{Z}(B^T)$ in $\mathbb{R}^k$. We may think of $\mathcal{P}(A)$ and $\mathcal{P}(B)$
as column-parallelotopes, and refer to $\mathcal{Z}(A^T)$ and
$\mathcal{Z}(B^T)$---which are defined using the columns of $A^T$ and
$B^T$---as the corresponding row-zonotopes (of $A$ and $B$). We ask for
conditions under which congruence of one pair of objects, coming from
equality of the corresponding shape matrices, implies congruence of the
other pair and how these conditions relate to the congruences. For example,

\medskip

{\footnotesize \lcrline{$\displaystyle$}{$A_1=\textstyle
\begin{bmatrix}
\hfill 5&\hfill 1\\
\hfill 1&\hfill 3
\end{bmatrix}
$, $B_1=\sqrt{2}
\begin{bmatrix}
\hfill 3&\hfill 2\\
\hfill 2&\hfill-1
\end{bmatrix}
$; $A_2=\frac{1}{13}
\begin{bmatrix}
\hfill 3&\hfill-12\\
\hfill 4&\hfill-3\\
\hfill 12&\hfill 4
\end{bmatrix}
$, $B_2=\frac{\sqrt{2}}{26}
\begin{bmatrix}
\hfill 15&\hfill-9\\
\hfill 7&\hfill 1\\
\hfill 8&\hfill 16
\end{bmatrix}
$;}{}

\medskip

\lcrline{} {$A_3=\textstyle
\begin{bmatrix}
\hfill 1&\hfill 2\\
\hfill 3&\hfill 4
\end{bmatrix}
$, $B_3=\frac{1}{\sqrt{884}}
\begin{bmatrix}
\hfill 46&\hfill 48\\
\hfill 82&\hfill 124
\end{bmatrix}
$; $A_4=
\begin{bmatrix}
\hfill 26&\hfill 8\\
\hfill 24&\hfill 2\\
\hfill 18&\hfill\,-16\\
\hfill 32&\hfill 26
\end{bmatrix}
$, $B_4=\sqrt{2}
\begin{bmatrix}
\hfill 17&\hfill 9\\
\hfill 13&\hfill 11\\
\hfill 1&\hfill 17\\
\hfill 29&\hfill 3
\end{bmatrix}
,$}{}}

\medskip

\noindent are four pairs of matrices where both $(A_i)^TA_i=(B_i)^TB_i$ and
$A_i(A_i)^T=B_i(B_i)^T$. In other words, $\mathcal{P}(A_i)$ is congruent to
$\mathcal{P}(B_i)$ and $\mathcal{Z}((A_i)^T)$ is congruent to
$\mathcal{Z}((B_i)^T)$.

Start by considering the case where $n=k$ and all four objects are
$n$-parallelotopes in $\mathbb{R}^n$. Consider the shape
matrix $A^TA$ of column-parallelotope $\mathcal{P}(A)$ and shape matrix
$AA^T$ of the corresponding row-parallelotope $\mathcal{P}(A^T)$. Another
matrix, $A^2=AA=(A^T)^TA$, can now be thought of as the {\bf comparison matrix}
between the generating vectors of the row-parallelotope and the
column-parallelotope. It seems plausible to conjecture that if matrices $A$
and $B$ have equal comparison matrices ($A^2=B^2$), the shape matrices of the
row-parallelotopes will be the same ($AA^T=BB^T$) if and only if the shape
matrices of the column-parallelotopes are the same ($A^TA=B^TB$). This is
in fact true.

\bigskip

\noindent {\bf Corollary 2.3.} {\it If $A$ and $B$ are square nonsingular matrices, and if
$A^2=B^2$, then $AA^T=BB^T$ if and only if $A^TA=B^TB$.}

\medskip

{\it Proof.} It suffices to prove one of the implications,
say$\Leftarrow\!.\,\,$Suppose $A^TA=B^TB$. From Proposition 2.1 it then
follows that $B=QA$, or $A=Q^TB$. Therefore, $AQ^TB=AA=BB$, which because
$B$ is nonsingular implies $AQ^T=B$, and so
$AA^T=AQ^T(AQ^T)^T=BB^T$.\hfill{\qedsymbol}

\medskip

Corollary 8.1 of [6] gives this same result in complex form
but makes no reference to the geometric interpretation involving row
and column-parallelotopes.

Another reasonable geometric conjecture is that if $A$ and $B$ have congruent
row-parallelotopes ($AA^T=BB^T$) and congruent column-parallelotopes
($A^TA=B^TB$), then their comparison matrices are identical ($A^2=B^2$) if
and only if the two congruences (provided by the matrix $Q$) are identical. This too is true.

\bigskip

\noindent {\bf Corollary 2.4.} {\it Let $A$ and $B$ be square nonsingular matrices such that
$AA^T=BB^T$ and $A^TA=B^TB$. Then $A^2=B^2$ iff there exists an
orthogonal matrix $Q$ such that $B=QA$ and $B^T=QA^T$.}

\medskip

{\it Proof.} Only $\Rightarrow$ is proved as $\Leftarrow \,$ is trivial. From
Proposition 2.1, $B=Q_1A$ and $B^T=Q_2A^T$. \ We must show that
$Q_1=Q_2.\,\,$Meanwhile, $A^2=B^2$ can be restated as $BA^{-1}=B^{-1}A$.
From the first and last of these $\,$several$\,$equalities,
$Q_1=BA^{-1}=B^{-1}A$. The second equality may also be rewritten as
$B=A(Q_2)^T=A(Q_2)^{-1}$, implying $Q_2=B^{-1}A$. Thus, $Q_1=B^{-1}A=Q_2$,
as was required. \hfill{\qedsymbol}

\medskip

When there is no orthogonal $Q$ such that $B=QA$ and $B^T=QA^T$ both hold, it
is possible to have $(1)$ $A^TA=B^TB$ and $(2)$ $AA^T=BB^T$, but $A^2\neq
B^2$. This happens, for example, in the case of the third pair of matrices
given above. In this situation the condition $A^2=B^2$ will be replaced with a weaker
comparison condition that does hold whenever $(1)$ and $(2)$ hold and
therefore seems more closely tied to these two conditions. Indeed, whenever
any two of $(1)$, $(2)$, and the new condition hold, it will turn out that
the third holds as well. Moreover, the condition will be defined and the
implications will hold in the more general setting of rectangular $n\times k$
matrices $A$ and $B$ with independent columns. In that case, the
column-parallelotopes $\mathcal{P}(A)$ and $\mathcal{P}(B)$ will reside in
$\mathbb{R}^n$ (with $n\geq k$) while the row-zonotopes $\mathcal{Z}(A^T)$
and $\mathcal{Z}(B^T)$ belong to $\mathbb{R}^k$. In order to obtain
comparison matrices in this setting, the parallelotopes $\mathcal{P}(A)$ and
$\mathcal{P}(B)$ are moved to congruent copies $\mathcal{P}(R)$ and $\mathcal{P}(S)$
within $\mathbb{R}^k$ by taking QR-decompositions $A=PR$ and $B=QS$
where $P$ and $Q$ are $n\times k$ matrices with orthonormal columns, and $R$ and $S$
are $k\times k$ upper triangular of rank $k$. By Theorem 2.2,
$\mathcal{P}(A)$ is indeed congruent to $\mathcal{P}(R)\subset\mathbb{R}^k$ and
$\mathcal{P}(B)$ is congruent to $\mathcal{P}(S)\subset\mathbb{R}^k$. The
parallelotopes in $\mathbb{R}^k$ can now be compared to the corresponding
row-zonotopes $\mathcal{Z}(A^T)$ and $\mathcal{Z}(B^T)$.

In order to make the comparison between $\mathcal{Z}(A^T)$ and $\mathcal{P}(R)$,
and between $\mathcal{Z}(B^T)$ and $\mathcal{P}(S)$, it does not suffice
to use $AR$ and $BS$. There are cases where $(1)$ and $(2)$ hold
but $AR=BS$ does not. An example is the pair $A_4$ and $B_4$ given above.
In order to make valid comparisons with the corresponding zonotopes, each
parallelotope must be allowed to independently reorient itself with respect to
its zonotope. For this purpose, additional orthogonal matrices $Q_1$ and $Q_2$
are introduced so that $\mathcal{Z}(A^T)$ is compared with
$\mathcal{P}(Q_1R)$ using the matrix $AQ_1R$, and $\mathcal{Z}(B^T)$ is
compared with $\mathcal{P}(Q_2S)$ using $BQ_2S$. Now everything works.
Setting

\medskip

\lcrline{}{$AQ_1R=BQ_2S$,}{$(3)$}

\medskip

\noindent it turns out that when any two of $(1)$, $(2)$, and the new condition $(3)$
hold, then so does the third. The precise relationship of $Q_1$ and $Q_2$ to $A,R,Q,$ and $S$
 will be clarified in the proof of Proposition 2.6, below.
 The following lemma will be used to establish the result.

\medskip

\noindent {\bf Lemma 2.5.} {\it Suppose $R$ and $S$ are non-singular $k\times k$ upper
triangular matrices with $R^TR=S^TS$. \ Then there is a diagonal matrix $J$ with $\pm 1$'s
on the diagonal such that $R=JS$.}

\medskip

{\it Proof.} \ We compute the first two rows of $R$ and $S$. A straightforward
induction (omitted) then completes the proof.

Let the columns of $R$ be $\boldsymbol{r}_1,\ldots,\boldsymbol{r}_k$ and
those of $S$ be $\boldsymbol{s}_1,\ldots,\boldsymbol{s}_k$. Let $m_{ij}$ be
the $ij^{\,th}$ entry of $R^TR=S^TS$. Then

\medskip

\lcrline{}{$(r_{11})^2=\boldsymbol{r}_1^T\boldsymbol{r}_1=m_{11}=%
\boldsymbol{s}_1^T\boldsymbol{s}_1=(s_{11})^2$,}{}

\medskip

\noindent from which $r_{11}=\pm s_{11}$. In addition, for each $j=2,\ldots,k$,

\medskip

\lcrline{}{$r_{11}r_{1j}=\boldsymbol{r}_1^T\boldsymbol{r}_j=m_{1j}=%
\boldsymbol{s}_1^T\boldsymbol{s}_j=s_{11}s_{1j}\,$.}{}

\medskip

\noindent If $r_{11}=s_{11}$, then $r_{1j}=s_{1j}$ for all $j=1,\ldots,k$ making the
first rows of $R$ and $S$ identical. If $r_{11}=-s_{11}$, then
$r_{1j}=-s_{1j}$ for all $j=1,\ldots,k$, so the first row of $R$ is the
negative of the first row of $S$.

Next, consider

\medskip

\lcrline{}{$(r_{12})^2+(r_{22})^2=\boldsymbol{r}_2^T\boldsymbol{r}_2=m_{22}=%
\boldsymbol{s}_2^T\boldsymbol{s}_2=(s_{12})^2+(s_{22})^2\,$.}{}

\medskip

\noindent We have already seen that $(r_{1j})^2=(s_{1j})^2$ for every $j$ including
$j=2$. It follows that$\,(r_{22})^2=(s_{22})^2$, or $r_{22}=\pm s_{22}$.
Meanwhile, for each $j=3,\ldots,k$,

\medskip

\lcrline{}{$\textunderset{${\scriptstyle =s_{12}s_{1j}}$}{$\underbrace{r_{12}r_{1j}}$}+r_{22}r_{2j} = %
\boldsymbol{r}_2^T\boldsymbol{r}_j = m_{2j} = %
\boldsymbol{s}_2^T\boldsymbol{s}_j=s_{12}s_{1j}+s_{22}s_{2j}\,$}{}

\smallskip

\noindent and so $r_{22}r_{2j}=s_{22}s_{2j}$. Consequently, either $r_{2j}=s_{2j}$
for every $j=2,\ldots,k$, or else $r_{2j}=-s_{2j}$ for every $j=2,\ldots,k$.
Therefore, the second rows of $R$ and $S$ are either identical or negatives
of each other. \ Continuing with similar computations, induction shows that
each row of $R$ is either the same or the negative of the corresponding row
of $S$. It follows that $R=JS$ as asserted.\hfill{\qedsymbol}

\medskip

From the lemma, if $A=QR=Q'R'$ are two QR-decompositions of an $n\times k$
matrix $A$ with independent columns, then $R'=JR$ and
$Q'=A(R')^{-1}=AR^{-1}J=QJ$. In other words, the QR-decomposition of $A$ is
unique up to a diagonal $k\times k$ matrix $J$ with $\pm 1$'s on the diagonal
(which changes the signs of specified columns of $Q$ and the corresponding
rows of $R$).

\medskip

\noindent {\bf Proposition 2.6.} {\it Let $A=PR$ and $B=QS$ be $n\times k$ matrices with
independent columns and QR-decompositions as indicated. \ Consider the three
conditions:

\smallskip

$(1)$ \ $A^TA=B^TB$,

\smallskip

$(2)$ \ $AA^T=BB^T$, and

\smallskip

$(3)$ \ there exist orthogonal matrices $Q_1$ and $Q_2$ such that
$\,AQ_1R=BQ_2S$.

\smallskip

If any two of the conditions hold, then so does the third. On the other
hand, no one of these conditions implies either of the other two.}

\medskip

{\it Proof.} $(a)$ Suppose $(3)$ and $A^TA=B^TB$ hold. It follows that
$R^TR=S^TS,$ so by Lemma 2.5, $R=JS$ (or $RS^{-1}=J$) for some diagonal matrix $J$
with $\pm 1$'s on the diagonal.  Condition $(3)$ then reduces to
$B=AQ_1JQ_2^T$, or $B^T=Q_2JQ_1^TA^T=Q^{\,\prime}A^T$ where
$Q^{\,\prime}=Q_2JQ_1^T$ is orthogonal, from which it follows that
$BB^T=AA^T$.

\smallskip

$(b)$ Suppose $(3)$ and $AA^T=BB^T$ hold. From $(3)$, it follows that
$B=AQ_1RS^{-1}Q_2^T$. Substituting the second equality in the first and
simplifying, $R^TR=S^TS$. Once again, Lemma 2.5 implies $S=JR$ for a diagonal
$J$ with $\pm 1$'s on the diagonal, so

\smallskip

\lcrline{}{$A^TA=R^TR=R^TJJR=S^TS=B^TB$.}{}

\smallskip

$(c)$ Finally, suppose conditions $(1)$ and $(2)$ hold. Then,
$A^TA=B^TB$ implies $R^TR=S^TS$, so by Proposition 2.1 there exists an
orthogonal $Q_1$ such that $S=Q_1R$. At the same time, applying Proposition
2.1 to $AA^T=BB^T$ guarantees existence of an orthogonal $Q_2$ such that
$B^T=Q_2A^T$. This last may be rewritten as $A=BQ_2$. It then follows that

\smallskip

\lcrline{}{$AQ_1R=AS=BQ_2S$.}{}

\smallskip

As for the last assertion of the proposition, it is clear that neither $(1)$
nor $(2)$ by itself implies either of the remaining two conditions. Giving
an example where $(3)$ holds but the other conditions do not will complete
the proof. To that end, let $A=
\begin{bmatrix}
\hfill 2&\hfill 6\\
\hfill 0&\hfill-1
\end{bmatrix}
$ and $B=
\begin{bmatrix}
\hfill 2&\hfill 2\\
\hfill 0&\hfill 1
\end{bmatrix}
$. \ Then $A^2=
\begin{bmatrix}
\hfill 4&\hfill 6\\
\hfill 0&\hfill 1
\end{bmatrix}
=B^2$, but $A^TA\neq B^TB$ and $AA^T\neq BB^T$. Meanwhile,
QR-decompositions for $A$ and $B$ may be taken as $A=PR=IA$ and $B=QS=IB$.
Also choosing $Q_1=\,Q_2=I$ then leads to $AQ_1R=A^2=B^2=BQ_2S$, so $(3)$
holds while $(1)$ and $(2)$ do not.\hfill{\qedsymbol}

\bigskip

We make several observations concerning the proposition. First, the
conclusion of the proposition may be rephrased as:

\medskip

{\it The pairs of conditions---$(1)\!+\!(3)$, $(2)\!+\!(3)$, and $(1)\!+\!(2)$---are
equivalent; but no one condition---$(1)$, $(2)$, or
$(3)$---implies either of the other two.}

\medskip

Second, if $A$ and $B$ are square matrices with QR-decompositions $A=PR$
and $B=QS$, and if condition (3) holds with $Q_1=P$ and $Q_2=Q$,
then condition $(3)$ becomes $A^2=B^2$.
When this version of the condition holds, Proposition 2.6 includes Corollary 2.3
and so is a generalization of that corollary.

Third, as $(3)$ must hold whenever $(1)$ and $(2)$ hold, the simultaneous occurrence
of congruences for both the column-parallelotopes and the row-zonotopes
ensures that $\,P(RQ_1R)=\,Q(SQ_2S)$.  This implies that the column spaces
of $P$ and $Q$ are the same and forces the column-parallelotopes
$\mathcal{P}(A)$ and $\mathcal{P}(B)$ to lie in the same $k$-dimensional
subspace of $\mathbb{R}^n$.

Fourth, denoting the rows of $A$ and $B\,$ as
$\boldsymbol{a}^1,\ldots,\boldsymbol{a}^n$ and
$\boldsymbol{b}^1,\ldots,\boldsymbol{b}^n$ respectively, and writing $R$ and
$S$ in terms of their columns as
$R=[\boldsymbol{r}_1,\ldots,\boldsymbol{r}_k]$ and
$S=[\boldsymbol{s}_1,\ldots,\boldsymbol{s}_k]$, condition $(3)$ becomes
$\boldsymbol{a}^i\cdot Q_1(\boldsymbol{r}_j)=\boldsymbol{b}^i\cdot
Q_2(\boldsymbol{s}_j)$ for every pair $(i,j)$. With $(2)$ and $(1)$ also
holding, $\left\|\boldsymbol{a}^i\right\|=\left\|\boldsymbol{b}^i\right\|$
and
$\left\|Q_1(\boldsymbol{r}_j)\right\|$=$\left\|Q_2(\boldsymbol{s}_j)\right%
\|$. Comparison condition $(3)$ then says that all corresponding angles
between the pairs $\left(\boldsymbol{a}^i,Q_1(\boldsymbol{r}_j)\right)$ and
$\left(\boldsymbol{b}^i,Q_2(\boldsymbol{s}_j)\right)$ are equal. (These are
the angles between the respective pairs of edges from the row-zonotope
$\mathcal{Z}(A^T)$ and reoriented column-parallelotope $\mathcal{P}(Q_1(R))$
on the one hand, and $\mathcal{Z}(B^T)$ and $\mathcal{P}(Q_2(S))$ on the
other.)

\medskip

QR-decompositions of $n\times k$ matrices with independent
columns, such as $A$ and $B$ with $A=PR$ and $B=QS$, produce ``generic''
parallelotopes $\mathcal{P}(R)$ and $\mathcal{P}(S)$ in $\mathbb{R}^k$,
independent of $n$. If the additional requirement is imposed on either
$R$ or $S$ that the entries on
its principal diagonal be positive, then that upper-triangular matrix is
uniquely determined. It represents a ``template'' parallelotope from which
all other congruent copies in $\mathbb{R}^n$ of that given shape of
parallelotope can be obtained by mapping using \text{$Q$-type} $n\times k$
matrices with $k$ orthonormal columns into appropriate $k$-dimensional
subspaces of Euclidean $n$-space $\mathbb{R}^n$.
A template parallelotope is a $k$-parallelotope in
$\mathbb{R}^k$ in ``standard position'' meaning that the \text{$j$-face}
defined by the first $j$ columns of the matrix---or by the $j$ corresponding
edges of the parallelotope---always lies in the subspace spanned by the first
$j$ standard basis vectors of the ambient space. The requirement that the
diagonal entries of the triangular matrix be positive implies, in addition,
that there exists a half-space such that the parallelotope and all of the
standard basis vectors $\boldsymbol{e}_1,\ldots,\boldsymbol{e}_n$ are
contained within that half-space.

\bigskip

\noindent {\bf Proposition 2.7.} {\it Suppose $A$ and $\,B$ are $n\times k$ matrices with
independent columns, $Q'$ is $m\times n$ with orthonormal columns, $A'=Q'A$,
and $B'=Q'B$. Then conditions $(1)$ and $(2)$ from Proposition 2.6 hold for
$A$ and $B$ if and only if the same conditions hold for $A'$ and $B'$.}

\medskip

{\it Proof.} $\Rightarrow\!$ : Suppose conditions $(1)$ and $(2)$ hold for $A$
and $B$. Then

\medskip

\lcrline{}{$A'\,^TA'=A^TQ'\,^TQ'A=A^TA=B^TB=B^TQ'\,^TQ'B=B'\,^TB'$}{}

\smallskip

\noindent and

\smallskip

\lcrline{}{$A'A'\,^T=Q'AA^TQ'\,^T=Q'BB^TQ'\,^T=B'B'\,^T$.}{}

\medskip

$\Leftarrow\!:$ \ \ \ Suppose conditions $(1)$ and $(2)$ hold for $A'$ and
$B'$. Then

\medskip

\lcrline{}{$A^TA=A^TQ'\,^TQ'A=A'\,^TA'=B'\,^TB'=B^TQ'\,^TQ'B=B^TB$,}{}

\noindent and

\smallskip

\lcrline{}{$AA^T=Q'\,^TA'A'\,^TQ'=Q'\,^TB'B'\,^TQ'=BB^T$.}{}  \hfill{\qedsymbol}

\bigskip

\section*{3. \,\,Volumes, Normal Vectors, and Rigidity of Zonotopes}

Symmetric cones, which were introduced in Section 1, can be used to derive a well-known volume formula for zonotopes in a new way.

\medskip

\noindent {\bf Proposition 3.1.}  {\it  Let $\mathcal{Z}(A)$ be an $n$-dimensional zonotope in
$\mathbb{R}^n$\,defined by an $n \times k$ matrix $A=[{\mathitbf a}_1,\ldots, {\mathitbf a} _k]$ of rank $n$
where the ${\mathitbf a}_j$'s are the columns of $A$. Then}

\medskip

\lcrline{} {$ \text{vol}_n \left(  \mathcal{Z} (A) \right) \,
= \,\, \displaystyle{\sum_{1 \leq j_1 < \cdots < j_n \leq k}} \,\,\, \left| \det\!\,\left( A^{j_1 , \ldots , \, j_n} \right) \right| \, $}{}

\smallskip

\noindent {\it where}  $A^{j_1 , \ldots , \, j_n}\,=\,[{\mathitbf a}_{j_1},\ldots, {\mathitbf a} _{j_n}]$.

\medskip

{\it Proof.} Central symmetry and convexity ensure that the zonotope decomposes completely into symmetric cones defined by pairs of opposite facets:

\medskip

\lcrline{} {$\mathcal{Z}(A)=
\displaystyle{\bigcup_{1\leq j_1<\cdots<j_{n-1}\leq k}}\{$cone$_{\boldsymbol{c}}(\mathcal{F}_{j_1,\ldots,\,j_{n-1}})\}$}{}

\smallskip

\noindent where $\boldsymbol{c}=\frac{1}{2}({\mathitbf a}_1+\cdots+{\mathitbf a}_k)$ is the center of symmetry of $\mathcal{Z}(A)$ and the facet
$\mathcal{F}_{j_1,\ldots,\,j_{n-1}}$ is one of a pair of translated copies of a generating facet defined by the $n \times (n-1)$ submatrix
$A^{j_1,\ldots,\,j_{n-1}}$. The generating facet is a zonotope of the form
$\mathcal{Z}(A^{j_1,\ldots,\,j_{n-1}})=\mathcal{Z}({\mathitbf a}_{j_1},\ldots, {\mathitbf a} _{j_{n-1}})$.
In degenerate cases, several such translated zonotopes might lie in the same hyperplane to form actual facets
of the given zonotope that are larger than parallelotopes, but this has no effect on the computation of volume.
For each submatrix of rank $n-1$, the normalized cross-product provides a unit normal vector for the corresponding
non-degenerate facet:

\smallskip

\lcrline{} {${\mathitbf n}_{j_1,\ldots,\,j_{n-1}}\,\,=\,\,
\displaystyle \frac{{\times}({\mathitbf a}_{j_1},\ldots, {\mathitbf a}_{j_{n-1}})}{|{\times}({\mathitbf a}_{j_1},\ldots, {\mathitbf a}_{j_{n-1}})|}
\,\,=\,\,\frac{{\times}({\mathitbf a}_{j_1},\ldots, {\mathitbf a}_{j_{n-1}})}{{\text {vol}}_{n-1}(\mathcal{F}_{j_1,\ldots,\,j_{n-1}})}.$}{}

\medskip

\noindent (Details about volumes defined by cross-products can be found, for example, in [4].)
The $n$-volume of each symmetric cone is $\frac{1}{n}$ times the $(n-1)$-volume of one of its antipodal bases times the height,
where the height is the distance between the pair of opposite bases of the cone.
That distance is simply the sum of the magnitudes of the projections of all $\mathitbf{a}_j$'s
onto ${\mathitbf n}_{j_1,\ldots,\,j_{n-1}}$. The magnitude of each projection is of the form

\medskip

\lcrline{} {$\displaystyle\left|\displaystyle {\mathitbf n}_{j_1,\ldots,\,j_{n-1}}\cdot\mathitbf{a}_j\right|
\,\,=\,\,\left|\frac{\times\left({\mathitbf a}_{j_1},\ldots, {\mathitbf a}_{j_{n-1}}\right)}{{\text {vol}}_{n-1}(\mathcal{F}_{j_1,\ldots,\,j_{n-1}})}\cdot\mathitbf{a}_j\right|
\,\,=\,\,\frac{\left| \det\!\,\left(A^{j_1,\ldots,\,j_{n-1},j}\right)\right|}{{\text {vol}}_{n-1}(\mathcal{F}_{j_1,\ldots,\,j_{n-1}})}$}{}

\medskip

\noindent where the second equality is obtained from the Laplace expansion of the determinant in its right-most column. The height is therefore

\medskip

\lcrline{} {$\displaystyle{\sum_{j\not= j_1,\ldots,j_{n-1}}}\,\,\,\frac{\left| \det \!\, (A^{j_1,\ldots,\,j_{n-1},j})\right|}{{\text {vol}}_{n-1}(\mathcal{F}_{j_1,\ldots,\,j_{n-1}})}
$}{}

\medskip

\noindent and the $n$-volume of the symmetric cone is

\begin{equation*}
\begin{aligned}\text{vol}_n\left(\text{cone}_{\boldsymbol{c}}(\mathcal{F}_{j_1,\ldots,\,j_{n-1}})\right)\,\,
   &= \,\,\displaystyle \frac{1}{n}\,{\text {vol}}_{n-1}(\mathcal{F}_{j_1,\ldots,\,j_{n-1}})\,\,\cdot\,\,{\sum_{j\not= j_1,\ldots,j_{n-1}}}\,\,\,\frac{\left|\text{\rm det}(A^{j_1,\ldots,\,j_{n-1},j})\right|}
{{\text {vol}}_{n-1}(\mathcal{F}_{j_1,\ldots,\,j_{n-1}})}\,\, \\[9pt]
   &= \,\,\frac{1}{n}\,\displaystyle{\sum_{j\not= j_1,\ldots,j_{n-1}}}\,\,\,\left| \det\!\,(A^{j_1,\ldots,\,j_{n-1},j})\right|.
\end{aligned}
\end{equation*}

\noindent It follows that

\begin{equation*}
\begin{aligned}\text{vol}_n\left(\mathcal{Z}(A)\right)\,\,
   &= \,\,\displaystyle{\sum_{1\leq j_1<\cdots<j_{n-1}\leq k}}\text{vol}_n\left(\text{cone}_{\boldsymbol{c}}(\mathcal{F}_{j_1,\ldots,\,j_{n-1}})\right) \\[9pt]
   &= \,\,\frac{1}{n}\,\sum_{\substack{1\leq j_1<\cdots<j_{n-1}\leq k \\ j\not= j_1,\ldots,j_{n-1}}}\,\,\,\left|\det\!\,\left(A^{j_1,\ldots,\,j_{n-1},j}\right)\right| \\[9pt]
   &= \,\,\sum_{1\leq j_1<\cdots<j_n\leq k}\,\,\,\left|\det\!\,\left(A^{j_1,\ldots,\,j_n}\right)\right|.
\end{aligned}
\end{equation*}

\smallskip

\noindent The last displayed equality holds because each term $\left|\det\!\,\left(A^{j_1,\ldots,\,j_{n-1},j}\right)\right|$ on the next-to-last line occurs
$n$ times in equivalent forms within that sum.  \hfill{\qedsymbol}

\bigskip

The volume formula also follows from the (non-unique) tiling of a zonotope into translated copies of its generating parallelotopes, each used exactly once.
This decomposition was cited in the proof of Corollary 1.7.

\bigskip

\noindent {\bf Proposition 3.2.} {\it Every $n$-dimensional zonotope formed from $k$ generating vectors in ${\mathbb R}^n$ decomposes into single translated copies
of each of its generating parallelotopes. These intersect each other only in lower-dimensional faces and together form a tiling of the zonotope by
$\binom{k}{n}$ parallelotopes.}

\smallskip

{\it Proof.} A one-dimensional zonotope (line segment) decomposes into subsegments that are translations of all of its generating line segments. And in all dimensions,
parallelotopes decompose trivially as themselves. Hence the proposition is true for all zonotopes of dimension $1$ and for zonotopes with $k=n$ generators
in any dimension $n$. Assume by induction that a decomposition of the required type exists for all zonotopes in every dimension $<n$
as well as for zonotopes with fewer than $k>n$ generators in dimenison $n$. Since $\mathcal{Z}(A)$ is $n$-dimensional, at least one of its subzonotopes generated by
$k-1$ column vectors is also $n$-dimensional. Thus, we may suppose without loss of generality that
$\mathcal{Z}(A)=\mathcal{Z}({\mathitbf a}_1,\ldots, {\mathitbf a} _k)=
\mathcal{Z}({\mathitbf a}_1,\ldots, {\mathitbf a} _{k-1})\,{\oplus}\,\mathcal{Z}({\mathitbf a}_k)$ where the
first summand is already $n$-dimensional.

Assume for now that none of the other generators lie in the 1-dimensional subspace spanned by ${\mathitbf a}_k$.
The visible surface of $\mathcal{Z}({\mathitbf a}_1,\ldots, {\mathitbf a} _{k-1})$ in the direction of ${\mathitbf a}_k$
consists (by Lemma 3.2 of [4]) of unique translates of all of the generating facets of the zonotope. Each facet is
defined by fewer than $k$ generators so by the induction assumption, each decomposes into unique copies of its
$(n-1)$-dimensional generating parallelotopes. Forming the Minkowski sum with $\mathcal{Z}({\mathitbf a}_k)$ has the effect of adding to
$\mathcal{Z}({\mathitbf a}_1,\ldots, {\mathitbf a} _{k-1})$ a zone of facets that all contain a translated copy of ${\mathitbf a}_k$, and a new
visible surface that is a copy of the old one
translated by $\mathitbf{a}_k$. It is always possible to fill the space between the original and translated copies of the visible surface
with $n$-parallelotopes whose bases are the $(n-1)$-dimensional parallelotopes from the decompositions of the facets of the visible surface
and whose remaining generating edge is, in every case, $\mathitbf{a}_k$.
Thus, in addition to the $\binom{k-1}{n}$ parallelotopes in the decomposition of
$\mathcal{Z}({\mathitbf a}_1,\ldots, {\mathitbf a} _{k-1})$, which exist by the induction assumption, there are
$\binom{k-1}{n-1}$ parallelotopes of the type just described, for a total of $\binom{k}{n}$ parallelotopes
that together form a decomposition of $\mathcal{Z}(A)$ of the required type.

In the case where several generators ${\mathitbf a}_{j+1},\ldots,{\mathitbf a}_k$ all lie in a single 1-dimensional subspace, convexity of
${\mathcal Z(A)}$ forces the edges defined by these generators to be contiguous. The sum of these generators then replaces ${\mathitbf a}_k$ in the
previous description. As a result, $(k-j)\cdot\binom{j}{n-1}$ distinct $n$-dimensional parallelotopes are created where the bases are
$(n-1)$-dimensional parallelotopes from the visible surface and the remaining generator is in turn
${\mathitbf a}_{j+1},\ldots,{\mathitbf a}_k$.

\hfill{\qedsymbol}

\bigskip

The  parallelotopes from the preceding proof that contain edge $\mathitbf{a}_k$ are sandwiched between the visible surface of
$\mathcal{Z}({\mathitbf a}_1,\ldots, {\mathitbf a} _{k-1})$ in the direction of ${\mathitbf a}_k$ and its translated copy
by $\mathitbf{a}_k$. Their union defines a partial shell of
parallelotopes forming an ``exterior wall'' of zonotope $\mathcal{Z}(A)=\mathcal{Z}({\mathitbf a}_1,\ldots, {\mathitbf a} _k)$, which Shephard [10] called
a {\bf cup of cubes}.
Labeling this as $\mathcal{C}_{1,\ldots\,k-1}(\mathitbf{a}_k )$, we have the decomposition
$\mathcal{Z}(A)\,=\,\mathcal{Z}({\mathitbf a}_1,\ldots, {\mathitbf a} _{k-1})\bigcup\mathcal{C}_{1,\ldots,k-1}(\mathitbf{a}_k )$. A similar
decomposition of the smaller zonotope yields

\medskip

\lcrline{}{$\mathcal{Z}(A)\,=\,\mathcal{Z}({\mathitbf a}_1,\ldots, {\mathitbf a}_{k-2})\bigcup\mathcal{C}_{1,\ldots,k-2}(\mathitbf{a}_{k-1})
\bigcup\mathcal{C}_{1,\ldots,k-1}(\mathitbf{a}_k ).$}{}

\noindent More generally,

\medskip

\lcrline{}{$\mathcal{Z}(A)\,=\,\mathcal{Z}({\mathitbf a}_{j_1},\ldots, {\mathitbf a}_{j_n})
\bigcup\mathcal{C}_{{j_1},\ldots,{j_n}}(\mathitbf{a}_{j_{n+1}})\bigcup\cdots\bigcup\mathcal{C}_{{j_1},\cdots,{j_{k-1}}}(\mathitbf{a}_{j_k} )$}{}

\medskip

\noindent where $\mathcal{Z}({\mathitbf a}_{j_1},\ldots, {\mathitbf a}_{j_n})$ is a generating parallelotope of $\mathcal{Z}(A)$
and $\mathcal{C}_{{j_1},\cdots,{j_l}}(\mathitbf{a}_{j_{l+1}})$ denotes the cup of cubes of $\mathcal{Z}({\mathitbf a}_{j_1},\ldots, {\mathitbf a}_{j_l})$
in $\mathcal{Z}({\mathitbf a}_{j_1},\ldots, {\mathitbf a}_{j_{l+1}})$ defined by $\mathitbf{a}_{j_{l+1}}$.
Thus, when $\mathcal{Z}(A)$ is developed as the Minkowski sum of successive line segments, various decompositions of the intermediate zonotopes
in the manner just shown produce different decompositions of $\mathcal{Z}(A)$ in terms of generating parallelotopes.
These decompositions, however, do not lead to all possible tilings of the zonotope. (See Shephard [10].)

\bigskip

To illustrate one possibility of what might happen to the facets and tiling of a zonotope, consider

\medskip

\[  A_0 \, = \, \left[\begin{array}{rrrrr}
      1 & 0 & 1 & 0 & -1 \\
      0 & 1 & 1 & 0 &  1 \\
      0 & 0 & 0 & 1 &  1
      \end{array} \right] \qquad \text{and}  \qquad   A_{\epsilon} \, = \,\left[\begin{array}{rrrrr}
                                   1 & 0 & 1 & 0 & -1 \\
                                   0 & 1 & 1 & 0 &  1 \\
                                   0 & 0 & \epsilon & 1 & 1
                                   \end{array}\right]. \]
\medskip

\noindent The first three columns of $A_0$ are dependent while the other triples of columns are independent. In $A_{\epsilon}$ for small non-zero ${\epsilon}$,
all ten triples of columns are independent. Consequently, only nine out of ten possible choices give generating parallelotopes of ${\mathcal Z}(A_0)$
while for ${\mathcal Z}(A_{\epsilon})$, all ten are parallelotopes. In each case, translates of the generating parallelotopes
can be arranged in various ways to form a
tiling of the zonotope. The generating subzonotope defined by the first three columns of $A_0$ is two-dimensional and translates to a pair of symmetrically
opposite hexagonal facets of the zonotope. Each of these facets coincides with a union of translations of the three generating facets
${\mathcal Z}({\mathitbf a}_1,{\mathitbf a}_2)$, ${\mathcal Z}({\mathitbf a}_1,{\mathitbf a}_3)$, and
${\mathcal Z}({\mathitbf a}_2,{\mathitbf a}_3)$, where ${\mathitbf a}_1,\ldots,{\mathitbf a}_5$ are the
columns of $A_0$. Zonotope ${\mathcal Z}(A_0)$ has $2\binom{5}{2}=20$ generating facets but only $14$ geometric facets; translates of three generating facets
make up each hexagonal facet. In the case of ${\mathcal Z}(A_{\epsilon})$, denote the columns of $A_{\epsilon}$ by ${\mathitbf a}'_1,\ldots,{\mathitbf a}'_5$.
Letting $\epsilon\!\rightarrow\!0$ demonstrates how a translate of the generating
parallelotope ${\mathcal Z}({\mathitbf a}'_1,{\mathitbf a}'_2,{\mathitbf a}'_3)$ flattens out and approaches one of the two
hexagonal facets on the boundary of  ${\mathcal Z}(A_0)$. Which facet is approached depends on the choice of tiling.
At the same time, two different sets of translates of the generating facets  ${\mathcal Z}({\mathitbf a}'_1,{\mathitbf a}'_2)$,
${\mathcal Z}({\mathitbf a}'_1,{\mathitbf a}'_3)$, and ${\mathcal Z}({\mathitbf a}'_2,{\mathitbf a}'_3)$ approach co-planarity
on opposite sides of the boundary of the zonotope to form copies of that same hexagonal facet.

\bigskip

Angles between edges $l\boldsymbol{a}_i$ and $l\boldsymbol{a}_j$ of an $n$-zonotope $\mathcal{Z}(A)\,=\,\mathcal{Z}(\boldsymbol{a}_1,\ldots,\boldsymbol{a}_k)$
can be computed as
$\theta=$arccos$(\boldsymbol{a}_i \cdot\boldsymbol{a}_j/ |\boldsymbol{a}_i|\,|\boldsymbol{a}_j|)$. Dihedral angles between facets can be computed
similarly using normal vectors to the facets. For facet
$\mathcal{Z}(A^{j_1,\ldots,j_{n-1}})\,=\,\mathcal{Z}(\boldsymbol{a}_{j_1},\ldots,\boldsymbol{a}_{j_{n-1}})$, the unit normal vector is

\medskip

\lcrline{} {${\mathitbf n}_{j_1,\ldots,\,j_{n-1}}\,\,=\,\,
\displaystyle \frac{{\times}({\mathitbf a}_{j_1},\ldots, {\mathitbf a}_{j_{n-1}})}{|{\times}({\mathitbf a}_{j_1},\ldots, {\mathitbf a}_{j_{n-1}})|}
\,\,=\,\,\frac{{\times}({\mathitbf a}_{j_1},\ldots, {\mathitbf a}_{j_{n-1}})}{{\text {vol}}_{n-1}(\mathcal{Z}(\boldsymbol{a}_{j_1},\ldots,\boldsymbol{a}_{j_{n-1}}))}.$}{}

\medskip

\noindent The cross-product used to find this normal vector can also be recovered from $\wedge^{n-1}(A)$, the matrix representing the map $\wedge^{n-1}f\!:\!\wedge^{n-1}\mathbb{R}^k\rightarrow\wedge^{n-1}\mathbb{R}^n$
with respect to  reverse lexicographically ordered bases of the exterior powers $\wedge^{n-1}\mathbb{R}^k$ and $\wedge^{n-1}\mathbb{R}^n$.
The entries of $\wedge^{n-1}(A)$ are the $(n-1) \times (n-1)$ minors of $A$, with the minor in row
$i$ and column $j$ the one defined by omitting those row and column indices. If $\wedge^{n-1}(A)$ is modified so that its rows and columns alternate in sign with
the $(i,j)$-th entry multiplied by $(-1)^{n+i+j}$, then the $j$-th column of the new matrix $\wedge^{n-1}_\pm(A)$ will be the cross product (or its negative) of the columns of $A$ with the complementary column indices. That is,
$$\big(\wedge^{n-1}_\pm(A)\big)^j=(-1)^{j}\big[{\times}({\mathitbf a}_1,\ldots,\widehat{\mathitbf a}_j ,\ldots ,{\mathitbf a}_n)\big].$$
The norm of the cross product then gives the $(n-1)$-volume of the parallelotope defined by those columns.
(See, for example, Corollary 1.3 of [4].
If the columns are not independent, they define a zonotope of rank less than $n-1$
whose $(n-1)$-volume is 0.) Thus, the columns of $\wedge^{n-1}_\pm(A)$ are normal vectors to the parallelotopes that comprise the generating facets of $\mathcal{Z}(A)$, and the norms of the vectors give the $(n-1)$-volumes of those parallelotopes. (Note that from Corollary 1.8, each facet decomposes into such parallelotopes.)

\bigskip

We wish to examine from the perspective of zonotopes two classic results in the theory of convex polytopes. The first is due to Minkowski. (See, for example, [5]):

\medskip

\noindent {\bf Theorem 3.3.} {\it \ Given distinct unit vectors
$\boldsymbol{u}_1 ,\ldots , \boldsymbol{u}_t$ that span $\mathbb{R}^n$ and corresponding arbitrary positive real numbers $a_1, \ldots , a_t$, then up to translation, there exists a unique convex polytope $\mathcal{P} \in \mathbb{R}^n$
for which the vectors are the outward-pointing normals to the facets and the numbers
are the $(n-1)$-volumes of the facets, if and only if \ $\sum a_i \boldsymbol{u}_i = \boldsymbol{0}$}.

\bigskip

The second, due to Cauchy in $\mathbb{R}^3$, extended to arbitrary $\mathbb{R}^n$ by Alexandrov (see, for instance, [1] ), and a basic part of geometric rigidity theory, is:

\bigskip

\noindent {\bf Theorem 3.4.} {\it \ If combinatorially equivalent convex polytopes in $\mathbb{R}^n$,
$n\geq 3$, have congruent corresponding facets, then the polytopes are congruent.}

\bigskip

Considering first Minkowski's theorem, observe that for polytopes whose facets come in pairs with equal $(n-1)$-volumes and unit normal vectors that are negatives of each other, the condition
$\sum a_i \boldsymbol{u}_i = \boldsymbol{0}$ is automatically satisfied. Indeed, it will turn out that given any distinct set of unit vectors spanning $\mathbb{R}^n$ and any corresponding set of positive reals, there exists a unique centrally-symmetric polytope whose pairs of opposite facets have the given unit vectors and their negatives as outward-pointing normals and the corresponding numbers as the common $(n-1)$-volumes of the pairs of facets. Thus, we have:

\bigskip

\noindent {\bf Proposition 3.5.} {\it Given distinct unit vectors
$\boldsymbol{u}_1 ,\ldots , \boldsymbol{u}_t$ spanning $\mathbb{R}^n$, no two of which are negatives of each other, and corresponding arbitrary positive real numbers $a_1, \ldots , a_t$, there exists a unique, centrally-symmetric polytope $\mathcal{P}$ with $2t$ facets such that $\boldsymbol{u}_i$ and $-\boldsymbol{u}_i$ are outward-pointing normals to facets $\mathcal{F}_i$ and
$\mathcal{F}^{\,\text{op}}_i$, and} vol$_{n-1}\mathcal{F}_i\, =\,\, $vol$_{n-1}\mathcal{F}^{\,\text{op}}_i =a_i$.

\medskip

{\it Proof.} The hypotheses ensure that the sum
$\sum a_i \boldsymbol{u}_i + \sum a_i (-\boldsymbol{u}_i) = \boldsymbol{0}$ when taken over the $(n-1)$-volumes of all facets times their corresponding normal vectors. Theorem 3.3 therefore applies and guarantees existence of a unique convex polytope $\mathcal{P}$ with the given vectors and their negatives as normal vectors and the given numbers as $(n-1)$-volumes of $t$ pairs of corresponding facets. The normal vectors come in opposite pairs, so opposite pairs of facets lie in parallel hyperplanes. The polytope is thus the intersection of the slabs that lie between these pairs of parallel hyperplanes. We need to show that
$\mathcal{P}$ is centrally symmetric.

In order to bound a closed polytope, there must be at least $n$ slabs. If there are exactly $n$ given vectors, they form a basis for $\mathbb{R}^n$, and
there are exactly $n$ slabs. $\mathcal{P}$ is then necessarily a parallelotope and therefore centrally symmetric. Thus, the proposition holds for $t=n$.

Assume, by induction, that the proposition holds for all sets of $<m$ normal vectors and corresponding facet volumes for a fixed value, $m$. Now consider $m$ unit normal vectors and $m$ corresponding
facet volumes $a_1, \ldots , a_m$. By Theorem 3.3, there exists a unique convex polytope $\mathcal{P}$ whose pairs of opposite facets satisfy the conditions of the proposition with these
values. We wish to show that this polytope is centrally symmetric.

$\mathcal{P}$ is the intersection of another polytope $\mathcal{P}\,'$ and a specific slab bounded by hyperplanes $\mathcal{H}_{m}$ and $\mathcal{H}_{-m}$,
which have outward-pointing normal vectors $\boldsymbol{u}_{m}$ and $\boldsymbol{u}_{-m}$. It is not clear, however, that $\mathcal{P}\,'$ satisfies the condition of the proposition requiring opposite
facets to have equal $(n-1)$-dimensional volumes. (In the end, it will turn out that $\mathcal{P}\,'$ does satisfy all of the conditions and is, in fact, centrally symmetric.)
To circumvent this difficulty, consider first an ``intermediate'' polytope $\mathcal{P}\,''$ that is the intersection of $\mathcal{P}\,'$
with the half-space defined by the hyperplane $\mathcal{H}_{-m}$ and its inward-pointing normal vector $-\boldsymbol{u}_{-m}$. In effect, $\mathcal{P}\,''$, which has one fewer facet than $\mathcal{P}$, is the polytope that results when the facet of $\mathcal{P}$ contained in hyperplane $\mathcal{H}_{m}$ is ``removed''. Compared with the remaining facets of $\mathcal{P}$, some of the corresponding facets
of $\mathcal{P}\,''$ have larger $(n-1)$-dimensional volumes, while the rest of the facets remain the same. Altering the list of values given for the facet volumes of $\mathcal{P}$ by eliminating the last value,
substituting the $(n-1)$-dimensional volumes of those facets from $\mathcal{P}\,''$ whose volumes increase compared to the corresponding facets of $\mathcal{P}$, and leaving the rest of the values unchanged,
a new list of numbers, $b_1, \ldots , b_{m-1}$, is obtained. By the inductional assumption, there is a unique centrally symmetric polytope $\mathcal{P}\,'''$ satisfying the conditions of the proposition with
specified unit normal vectors $\boldsymbol{u}_1 ,\ldots , \boldsymbol{u}_{m-1}$ and their opposites, and with corresponding pairs of facet volumes $b_1, \ldots , b_{m-1}$. Denote the center of
$\mathcal{P}\,'''$ by $\boldsymbol{c}$.

$\mathcal{P}$ will equal the intersection of $\mathcal{P}\,'''$ with a particular slab bounded by hyperplanes $\mathcal{H}_{m}$ and $\mathcal{H}_{-m}$
whose outward-pointing normal vectors are $\boldsymbol{u}_{m}$ and $\boldsymbol{u}_{-m}$. All slabs are centally symmetric. If it can be demonstrated that this slab has $\boldsymbol{c}$ as a center,
then $\mathcal{P}$ will be the intersection of two centrally symmetric sets with the same center and therefore will be centrally symmetric with respect to $\boldsymbol{c}$ by Lemma 1.1$(f)$.
To see that the slab bounded by $\mathcal{H}_{m}$ and $\mathcal{H}_{-m}$ is centrally symmetric with respect to $\boldsymbol{c}$, let $\mathcal{H}$ be
the supporting hyperplane of $\mathcal{P}\,'''$ for which $\boldsymbol{u}_m$ is the outward-pointing normal, and let $\mathcal{H}\,'$ be the parallel supporting hyperplane on the opposite side
of $\mathcal{P}\,'''$. The distance $w$ between the hyperplanes is the width of $\mathcal{P}\,'''$ in the
direction of $\boldsymbol{u}_m$. Let $\mathcal{H}_{tw}$ be the hyperplane parallel to and between $\mathcal{H}$
and $\mathcal{H}\,'$ whose distance from $\mathcal{H}$ is $tw$, where $0\leq t\leq 1$. Define a non-negative
real-valued function $C\!:\![0,1]\rightarrow \mathbb{R}$ such that $C(t)$ is the $(n-1)$-dimensional cross-sectional
volume of $\mathcal{P}\,'''\cap \mathcal{H}_{tw}$. As $\mathcal{P}\,'''$ is centrally symmetric,
this function is unimodal and symmetric about the value $t=\frac{1}{2}$, where it attains its
maximum. Moreover, the cross section at $t=\frac{1}{2}$ contains the center $\boldsymbol{c}$ of
$\mathcal{P}\,'''$. (All of these follow from Corollary 2.2 of [2].) As a consequence, every slab bounded by hyperplanes of the form $\mathcal{H}_t$ and $\mathcal{H}_{1-t}$ will be centrally symmetric
with respect to $\boldsymbol{c}$, and the intersection of each such slab with $\mathcal{P}\,'''$ will be a centrally symmetric polytope.
In particular, this will be the case for the value $t=t_0$ where the $(n-1)$-dimensional cross-sectional volumes of $\mathcal{P}\,'''\cap \mathcal{H}_{t_0}$ and
$\mathcal{P}\,'''\cap \mathcal{H}_{1-{t_0}}$ are both $a_m$, and where $\mathcal{H}_{t_0}=\mathcal{H}_m$ and $\mathcal{H}_{1-t_0}=\mathcal{H}_{-m}$.
The values obtained for the $(n-1)$-dimensional volumes of facets of the intersection polytope that have non-empty intersections with $\mathcal{H}_{t_0}$ must agree with the corresponding facet volumes
of $\mathcal{P}$ because these facets are formed exactly as facets of $\mathcal{P}$ had been formed by intersecting $\mathcal{P}\,'$ with a similar slab. It follows that all facet volumes agree with those of
$\mathcal{P}$, and hence the intersection of $\mathcal{P}\,'''$ with the slab between $\mathcal{H}_{t_0}$ and $\mathcal{H}_{1-t_0}$
is in fact $\mathcal{P}$, which is therefore seen to be centrally symmetric.
(This also shows that $\mathcal{P}\,'=\mathcal{P}\,'''$ and that the slab whose intersection with $\mathcal{P}\,'$ produced $\mathcal{P}$ is the same as the slab between $\mathcal{H}_{t_0}$
and $\mathcal{H}_{1-t_0}$\,).
\hfill{\qedsymbol}

\bigskip

\noindent {\bf Corollary 3.6.}{\it \ If $\{\boldsymbol{u}_1 ,\ldots , \boldsymbol{u}_n\}$ is a basis for
$\mathbb{R}^n$ and $a_1,\ldots ,a_n$ are arbitrary positive real numbers, then there exists a unique parallelotope with pairs of opposite facets having the $\boldsymbol{u}_i$'s as normal vectors and the $a_i$'s as the $(n-1)$-volumes of the facets.}

\medskip

{\it Proof.} The first step in the proof of the proposition included the observation that when the vectors $\{\boldsymbol{u}_1 ,\ldots , \boldsymbol{u}_n\}$ were a basis for $\mathbb{R}^n$, the polytope uniquely
determined by Minkowski's Theorem using the arbitrarily given $a_1,\ldots ,a_n$ is bounded by $n$ slabs and is necessarily a parallelotope. \hfill{\qedsymbol}

\bigskip

\noindent As a result, parallelotopes can also be found with arbitrary dihedral angles $0<\theta <2\pi$ and facet volumes.
In the simple case of boxes, one might give high school students taking elementary algebra the problem of finding the dimensions of a box whose opposite pairs of faces have areas 1, 2, and 3, or, for that
matter, any three arbitrarily picked positive areas.

We also note the following algebraic consequence of the preceding geometric corollary:

\bigskip

\noindent {\bf Corollary 3.7.}{\it \ Every $n\times n$ non-singular real matrix, $B$, is the $(n-1)$-st exterior power of a unique $n\times n$ matrix. That is, $B=\wedge^{n-1}(A)$ for some $n\times n$ matrix $A$, which may be regarded as the $(n-1)$-st exterior root of $B$.}

\medskip

{\it Proof.} Let $B$ be an arbitrary non-singular $n\times n$ matrix. Regard the columns of $B$ and their negatives as the normal vectors to pairs of facets of a parallelotope where the norms of these vectors are the corresponding facet-volumes. Recall from the discussion preceding Theorem 3.3 that for an $n\times n$ matrix $A$, the $j$-th column of $\wedge^{n-1}_\pm(A)$ is
$(-1)^j\big[{\times}({\mathitbf a}_1,\ldots,\widehat{\mathitbf a}_j ,\ldots ,{\mathitbf a}_n)\big]$. This column vector is the outward normal for one of the pair of facets of the parallelotope $\mathcal{P}(A)$ whose corresponding generating facet is
$\mathcal{P}({\mathitbf a}_1,\ldots,\widehat{\mathitbf a}_j ,\ldots ,{\mathitbf a}_n)$.
The norm of the vector is the $(n-1)$-volume of this facet. The previous corollary guarantees existence of a unique parallelotope $\mathcal{P}(A)$ with specified outward-pointing normals and facet-volumes. Its defining matrix $A$, with perhaps the columns permuted and some of their signs changed, then satisfies $\wedge^{n-1}_\pm(A)=B$. If $B$ is first altered to $B'$ where the $(i,j)$-th entry of $B'$ is $(-1)^{n+i+j}$ times the corresponding entry of $B$, then $\wedge^{n-1}(A)=B$.  \hfill{\qedsymbol}

\bigskip

\noindent Conditions on matrices of various shapes that guarantee the existence of $k$-th exterior roots for different values of $k$ are less known.

\smallskip

We now consider the Cauchy-Alexandrov Theorem for zonotopes, where a direct proof is possible.

\bigskip

\noindent {\bf Proposition 3.8.}{\it \ If combinatorially equivalent zonotopes in $\mathbb{R}^n$, $n\geq 3$, have congruent corresponding facets, then the zonotopes are congruent.}

\medskip

{\it Proof.} Consider combinatorially equivalent zonotopes, $\mathcal{Z}(A')$ and $\mathcal{Z}(B')$, with respective defining matrices $A'$ of shape $p\times k$ and $B'$ of shape $q\times k'$.
Suppose $p\leq q$. After replacing $A'$ with $QA'$ where $Q$ is a $q\times p$ matrix with orthonormal columns, we may suppose, without loss of generality, that the zonotopes
are embedded in the same Euclidean space $\mathbb{R}^q$. Combinatorial equivalence means the zonotopes have the same
face lattice structure, and hence both have facets of the same dimension (and therefore defining matrices of the same rank $n\leq p$), as well as
the same number of edges (so that $k=k'$). As both are embedded in $n$-dimensional subspaces of $\mathbb{R}^q$, we may assume that both reside within $\mathbb{R}^n$
and are defined by matrices $A$ and $B$ of rank $n$ and shape $n\times k$. Thus, it suffices to consider combinatorially equivalent $n$-zonotopes $\mathcal{Z}(A)$ and $\mathcal{Z}(B)$
with congruent corresponding facets contained in $\mathbb{R}^n$ whose respective defining matrices, $A$ and $B$, are both of rank $n$ and shape $n\times k$. We wish to show that these zonotopes are
congruent to each other.

Consider first  the case where $k=n$. The columns of both matrices are then independent and the zonotopes  they define are parallelotopes. Corresponding generating facets of the parallelotopes are of the form
$\mathcal{P}({\mathitbf a}_1,\ldots,\widehat{\mathitbf a}_j ,\ldots ,{\mathitbf a}_n)$ and
$\mathcal{P}({\mathitbf b}_1,\ldots,\widehat{\mathitbf b}_j ,\ldots ,{\mathitbf b}_n)$, which can be
denoted briefly as $\mathcal{P}(A(\widehat{\jmath}\,))$ and $\mathcal{P}(B(\widehat{\jmath}\,))$.
As the facets in each pair are congruent, Theorem 2.2 implies that
$(A(\widehat{\jmath}\,))^T(A(\widehat{\jmath}\,))=(B(\widehat{\jmath}\,))^T(B(\widehat{\jmath}\,)).$
Taken over all values $j=1,\ldots n$, all columns of the matrices are covered by equivalent comparisons, and so $A^TA=B^TB$. (Indeed, it suffices to consider just three congruent pairs of corresponding facets in order to guarantee that all corresponding pairs of edges of the parallelotopes have been compared. This will be made explicit right after the end of the proof.)
By Theorem 2.2, the parallelotopes are therefore congruent to each other.

The proof will be completed by induction. Assume the proposition holds for all pairs of $n$-zonotopes with congruent corresponding facets and defining matrices of shape $n\times j$ where $j < k$.
Consider $n$-zonotopes $\mathcal{Z}(A)$ and $\mathcal{Z}(B)$ with congruent corresponding facets and defining matrices of shape $n\times k$.
Let $A(\widehat{k}\,)$ and $B(\widehat{k}\,)$ be the corresponding matrices with $k$-th columns omitted. It follows that
$\mathcal{Z}(A)=\mathcal{Z}(A(\widehat{k}\,))\,\oplus\, l\boldsymbol{a}_k$ and $\mathcal{Z}(B)=\mathcal{Z}(B(\widehat{k}\,))\,\oplus\, l\boldsymbol{b}_k$
where $l\boldsymbol{a}_k$ and $l\boldsymbol{b}_k$ are the edges defined by the last columns.
By the inductional assumption, $\mathcal{Z}(A(\widehat{k}\,))$ and $\mathcal{Z}(B(\widehat{k}\,))$ are congruent to each other, and as congruence extends downward to lower-dimensional faces,
$l\boldsymbol{a}_k$ and $l\boldsymbol{b}_k$ are of equal length.
Congruence of $\mathcal{Z}(A(\widehat{k}\,))$ and $\mathcal{Z}(B(\widehat{k}\,))$ implies that their visible surfaces in the respective directions
defined by $\boldsymbol{a}_k$ and $\boldsymbol{b}_k$ are also congruent. Moreover, the argument given in the preceding paragraph guarantees congruence for every pair of corresponding parallelotopes
that tile the zonotopes $\mathcal{Z}(A)$ and $\mathcal{Z}(B)$. Consequently, the cups of cubes,
$\mathcal{C}_{1,\ldots ,k-1}(\boldsymbol{a}_k)$ and $\mathcal{C}_{1,\ldots ,k-1}(\boldsymbol{b}_k)$ (defined in the comments following the proof of Proposition 3.2),
are also congruent. As $\mathcal{Z}(A)=\mathcal{Z}(A(\widehat{k}\,))\,\cup\, \mathcal{C}_{1,\ldots ,k-1}(\boldsymbol{a}_k)$
and $\mathcal{Z}(B)=\mathcal{Z}(B(\widehat{k}\,))\,\cup\, \mathcal{C}_{1,\ldots ,k-1}(\boldsymbol{b}_k)$, it follows that $\mathcal{Z}(A)$ and $\mathcal{Z}(B)$
are congruent to each other. \hfill{\qedsymbol

\bigskip

In the second paragraph of the proof, it was asserted that when three pairs of corresponding facets of combinatorially equivalent parallelotopes are congruent,
then the parallelotopes are themselves congruent. In terms of comparisons of the defining matrices, this is the same as saying that
if $(A(\widehat{\jmath}\,))^T(A(\widehat{\jmath}\,))=(B(\widehat{\jmath}\,))^T(B(\widehat{\jmath}\,))$ holds for three different index values, then $A^TA=B^TB$.
To see this in more detail, observe that $A^TA=B^TB$ is equivalent to $({\mathitbf a}_i)^T({\mathitbf a}_j)=({\mathitbf b}_i)^T({\mathitbf b}_j)$ holding for every pair $(i,j)$ with $1\leq i,j\leq n$.
If a pair of corresponding facets from $\mathcal{P}(A)$ and $\mathcal{P}(B)$ omit the edge defined by column $i_0$ of matrices $A$ and $B$ respectively,
then congruence of these facets is equivalent to $(A(\widehat{i_0}\,))^T(A(\widehat{{i_0}\,}\,))=(B(\widehat{i_0}\,))^T(B(\widehat{i_0}\,))$,
which in turn says that the comparisons $({\mathitbf a}_i)^T({\mathitbf a}_j)=({\mathitbf b}_i)^T({\mathitbf b}_j)$ hold for all $i$ and $j$ except for
$({\mathitbf a}_{i_0})^T({\mathitbf a}_j)=({\mathitbf b}_{i_0})^T({\mathitbf b}_j)$ and $({\mathitbf a}_i)^T({\mathitbf a}_{i_0})=({\mathitbf b}_i)^T({\mathitbf b}_{i_0})$.
If a second pair of corresponding facets are congruent and omit the edge defined by column $j_0$, then the only omitted comparisons this time are those of the form
$({\mathitbf a}_{j_0})^T({\mathitbf a}_j)=({\mathitbf b}_{j_0})^T({\mathitbf b}_j)$ and $({\mathitbf a}_i)^T({\mathitbf a}_{j_0})=({\mathitbf b}_i)^T({\mathitbf b}_{j_0})$.
If both pairs of facets are congruent, then taken together, the only missing comparisons are $({\mathitbf a}_{i_0})^T({\mathitbf a}_{j_0})=({\mathitbf b}_{i_0})^T({\mathitbf b}_{j_0})$
and $({\mathitbf a}_{j_0})^T({\mathitbf a}_{i_0})=({\mathitbf b}_{j_0})^T({\mathitbf b}_{i_0})$. Now, if a third pair of corresponding facets are congruent and omit a third edge defined
by column $k_0$, then while some other comparisons might be missing from just this congruence, when taken together with the other congruences, the missing comparisons
from the first two are now included as part of of the third congruence. As a result, all comparisons needed to establish $A^TA=B^TB$ are made when three different pairs of
corresponding facets from the parallelotopes are congruent, so in that case, the parallelotopes will be congruent to each other.

It becomes less obvious and perhaps more surprising that congruence of three pairs of facets still suffices in dimensions much greater than 3. It should also be noted that requiring the dimension to be at least 3 is necessary in order to guarantee that comparison of shape matrices for three corresponding pairs of facets is sufficient to imply $A^TA=B^TB$. When $n=2$,
$$({\mathitbf a}_1)^T({\mathitbf a}_1)=({\mathitbf b}_1)^T({\mathitbf b}_1) \text{\ \, and \ \,}
({\mathitbf a}_2)^T({\mathitbf a}_2)=({\mathitbf b}_2)^T({\mathitbf b}_2)$$ do not
force $$[{\mathitbf a}_1,{\mathitbf a}_2]^T[{\mathitbf a}_1,{\mathitbf a}_2]=[{\mathitbf b}_1,{\mathitbf b}_2]^T[{\mathitbf b}_1,{\mathitbf b}_2]$$
because they lack the comparison
$({\mathitbf a}_1)^T({\mathitbf a}_2)=({\mathitbf b}_1)^T({\mathitbf b}_2)$. Indeed, zonogons in $\mathbb{R}^2$ with congruent corresponding edges need not be congruent.

\medskip

For combinatorially equivalent zonotopes, the uniqueness part of Minkowski's Theorem can be proven directly. Theorem 3.3 and its consequences will not be used, but Corollary 3.7, which has an independent algebraic proof, will be.

\medskip

\noindent {\bf Proposition 3.9.}{\it \ If combinatorially equivalent zonotopes $\mathcal{Z}(A)$ and $\mathcal{Z}(B)$ have corresponding equal unit facet normals and facet volumes, then they are congruent.}

\medskip

{\it Proof.} Let $A=[\boldsymbol{a}_1 ,\ldots , \boldsymbol{a}_k]$ and $B=[\boldsymbol{b}_1 ,\ldots , \boldsymbol{b}_k]$ be matrices of rank $n$ with $k\geq n$, and suppose
these matrices define combinatorially equivalent zonotopes $\mathcal{Z}(A)$ and $\mathcal{Z}(B)$ in $\mathbb{R}^n$.
The number of generating facets of each zonotope will be the number of maximal subsets of columns of rank $(n-1)$. That is, a generating facet will be defined by a subset of
columns of rank $(n-1)$ to which no further columns can be added without increasing the rank. (Each generating facet will produce two bounding facets of a zonotope.)
The number of bounding facets will thus be some number $2m$ where $n\leq m\leq \binom{k}{n-1}$.
We will be less interested in the bounding facets themselves than in the parallelotope constituents of those facets defined by choosing exactly $(n-1)$ corresponding columns
from each defining matrix. While $m$ unit vectors and their negatives will represent the outward-pointing normals of the bounding facets for each zonotope,
repeating a normal vector for every parallelotope constituent of a bounding facet will produce a total number of $t\,{\mathrel{\mathop :}=}\,\binom{k}{n-1}$
normal vectors, $\boldsymbol{u}_1 ,\ldots , \boldsymbol{u}_t$, along with their negatives, that will be used in the description of each zonotope.
An equal number of non-negative numbers, $a_1,\ldots , a_t$, will represent the $(n-1)$-volumes of the corresponding pairs of parallelotope facet-constituents of each zonotope.

The proof is by induction on $k$. When $k=n$, the matrices are non-singular and define parallelotopes $\mathcal{P}(A)$ and $\mathcal{P}(B)$.
Corresponding facets of the parallelotopes have the same normal vectors, so the dihedral angles between pairs of facets are also equal. Corresponding facets also have the same volumes. Taken together, these comparisons ensure that the defining matrices satisfy
$$(\wedge^{n-1}_\pm A)^T(\wedge^{n-1}_\pm A)=(\wedge^{n-1}_\pm B)^T(\wedge^{n-1}_\pm B),$$
and therefore also  $$(\wedge^{n-1}A)^T(\wedge^{n-1}A)=(\wedge^{n-1}B)^T(\wedge^{n-1}B).$$

\smallskip

\noindent Moreover, $(\wedge^{n-1}A)^T(\wedge^{n-1}A)=\wedge^{n-1}(A^{T}A)$, from which
$$\wedge^{n-1}(A^{T}A)=\wedge^{n-1}(B^{\,T}B).$$
Corollary 3.7 then implies $$A^{T}A=B^TB,$$
so by Theorem 2.2, $\mathcal{P}(A)$ and $\mathcal{P}(B)$ are congruent.

Now assume the proposition holds for zonotopes defined by matrices with fewer than $k$ columns
for some fixed value $k>n$. Suppose $\mathcal{Z}(A)$ and $\mathcal{Z}(B)$ satisfy the normal-vector and facet-volume conditions
and are defined by $n\times k$ matrices. Let $A(\widehat{k}\,)$ and $B(\widehat{k}\,)$ be the corresponding
matrices with $k$-th columns omitted. It follows that
$\mathcal{Z}(A)=\mathcal{Z}(A(\widehat{k}\,))\,\oplus\, l\boldsymbol{a}_k$ and
$\mathcal{Z}(B)=\mathcal{Z}(B(\widehat{k}\,))\,\oplus\, l\boldsymbol{b}_k$. The normal vector to each facet of $\mathcal{Z}(A)$ belonging to the zone (that is, 1-zone) of facets containing
$\boldsymbol{a}_k$ is orthogonal to $\boldsymbol{a}_k$. All of these vectors span a hyperplane
orthogonal to $\boldsymbol{a}_k$. A similar relationship holds in $\mathcal{Z}(B)$. As the normal vectors and hyperplanes are the same, it follows that
$\boldsymbol{a}_k$ and $\boldsymbol{b}_k$ are parallel. Moreover, corresponding facets in the zones for $\boldsymbol{a}_k$ and $\boldsymbol{b}_k$ have the same volumes. The facets in these zones are
Minkowski sums of faces from either $A(\widehat{k}\,)$ with $l\boldsymbol{a}_k$ or from
$B(\widehat{k}\,)$ with $l\boldsymbol{b}_k$\,, respectively. The faces are either $(n-2)$- or $(n-1)$-dimensional,
and the resulting facets, after forming the sums, are then either prisms in the first case, or convex hulls of translated facets (when $\boldsymbol{a}_k$ or
$\boldsymbol{b}_k$ lies in the hyperplane containing the facet) in the second case. In either case,
congruence of the corresponding base faces, the fact that $\boldsymbol{a}_k$ and $\boldsymbol{b}_k$ are parallel, and equality of volumes of
the resulting facets, force $\boldsymbol{a}_k$ and $\boldsymbol{b}_k$ to have the same length. Once the
vectors are parallel and of the same length, the corresponding facets formed as Minkowski sums using these
vectors are congruent. Thus, all corresponding pairs of facets from $\mathcal{Z}(A)$ and $\mathcal{Z}(B)$
are congruent, and the two zonotopes are themselves congruent by Proposition 3.8.  \hfill{\qedsymbol}

\bigskip

\vspace{3.5in}

\bigskip

\noindent {\bf REFERENCES}

\medskip

\footnotesize

\begin{tabbing}

\hspace{0.35in}\=\kill

\lbrack 1\rbrack\LTab{A. D. Alexandrov,
{\it Convex Polyhedra}, Springer-Verlag, New York, 2005.}\\[5pt]

\lbrack 2\rbrack\LTab{D. Avis, et. al.,
{\it On the sectional area of convex polytopes}, Proceedings of the XII Annual Symposium\\
\LTab on Computational Geometry, New York, 1996.}\\[5pt]

\lbrack 3\rbrack\LTab{A. Bj\"orner, et. al.,
{\it Oriented Matroids} (2nd ed.), Cambridge U. Press, New York, 1999.}\\[5pt]

\lbrack 4\rbrack\LTab{E. Gover and N. Krikorian, {\it Determinants and the volumes of
parallelotopes and zonotopes}},\\
\LTab{Linear Algebra and its Applications, $\,\boldsymbol{433}\,$(2010), 28--40.}\\[5pt]

\lbrack 5\rbrack\LTab{B. Gr\"unbaum,
{\it Convex Polytopes}, Springer-Verlag, New York, 2003.}\\[5pt]

\lbrack 6\rbrack
\LTab{R. Horn and I. Olkin, {\it When does $A^*\!A\,=\,B^*\!B$
and why does one want to know?}}\\
\LTab{MAA Monthly, $\,\boldsymbol{103}\,$(1996), 470--482.}\\[5pt]

\lbrack 7\rbrack\LTab{P. McMullen, {\it Polytopes with centrally symmetric
faces},$\,$Israel$\,$J.$\,$Math.$\,\boldsymbol{8}\,$(1970), 194--196.}\\[5pt]

\lbrack 8\rbrack\LTab{P. McMullen, {\it Polytopes with centrally symmetric
facets},$\,$Israel$\,$J.$\,$Math.$\,\boldsymbol{23}\,$(1976), 337--338.}\\[5pt]

\lbrack 9\rbrack\LTab{G. Shephard, {\it Polytopes with centrally symmetric
faces},$\,$Canadian$\,$J.$\,$Math.$\,\boldsymbol{19}\,$(1967),
1206--1213.}\\[5pt]

\lbrack 10\rbrack\LTab{G. Shephard,
{\it Combinatorial$\,$properties$\,$of$\,$associated$\,$zonotopes},$\,$Canadian$\,%
$J.$\,$Math.$\,\boldsymbol{26}\,$(1974), 302--321.}

\end{tabbing}

\vspace{0.75in}

\normalsize

\begin{tabbing}

\hspace{0.35in}\=\kill

\LTab{Eugene Gover}\\

\LTab{Department of Mathematics}\\

\LTab{Northeastern University}\\

\LTab{Boston, MA 02115, U.S.A.}\\

\LTab{e.gover@neu.edu}

\end{tabbing}

\end{document}